\documentstyle[12pt]{article}

\newcommand{\qed}{\hfill \mbox{$\Box$ \ } \bigskip}

\newtheorem{lemma}{{\sc Lemma}}[section]
\newtheorem{thm}{{\sc Theorem}}[section]
\newtheorem{proposition}{{\sc Proposition}}[section]

\makeatletter 
 
\@addtoreset{equation}{section}
\makeatother

\newcommand{\bysame}{%
  \leavevmode\hbox to 3em{\hrulefill}\,}

\begin{document}

\title{Asymptotic expansion of the Bergman kernel \\
for weakly pseudoconvex tube domains in {\bf C}$^2$}
\author{by \\ Joe {\sc Kamimoto} \\
{\small Graduate School of Mathematical Sciences,  
The University of Tokyo,} \\
{\small 3-8-1, Komaba, Meguro, 
Tokyo, 153 Japan.}\\
{\it E-mail} : {\tt kamimoto@ms.u-tokyo.ac.jp} }
\date{October 10, 1996} 
\maketitle

%\begin{large} 

%\begin{center}\rule[.3cm]{14cm}{0.2mm}\end{center}

%\vspace{1 em}

\begin{abstract}
\footnote{{\it Math Subject Classification.} 32A40, 32F15, 32H10.} 

In this paper we give an asymptotic expansion 
of the Bergman kernel for certain 
weakly pseudoconvex tube domains of finite type in ${\bf C}^2$. 
Our asymptotic formula asserts that  
the singularity of the Bergman kernel 
at weakly pseudoconvex points 
is essentially expressed by using two variables ;  
moreover certain real blowing-up is 
necessary to understand its singularity. 
The form of the asymptotic expansion 
with respect to each variable  
is similar to that 
in the strictly pseudoconvex case 
due to C. Fefferman. 
We also give an analogous result in the case of the 
Szeg\"o kernel. 
\end{abstract} 
%\vspace{3em}

%%%%%%%%%%%%%%%%%%%%%%%%%%%%%%%%%%%%%%%%%%%%%%%%%%%%%%%%%

\section{Introduction}

The purpose of this paper is to 
give an asymptotic expansion 
of the Bergman kernel for certain class of  
weakly pseudoconvex tube domains of finite type in ${\bf C}^2$. 
We also give an analogous result of the Szeg\"o kernel for the 
same class of tube domains. 

%\vspace{.5 em}

Let $\Omega$ be a  domain with smooth boundary 
in ${\bf C}^n$. 
The Bergman space $B(\Omega)$ is the subspace of $L^2(\Omega)$ 
consisting of 
holomorphic $L^2$-functions on $\Omega$. 
The Bergman projection is the orthogonal projection 
${\bf B}:L^2(\Omega)\to B(\Omega)$. 
We can write ${\bf B}$ as an integral operator 
$$
{\bf B}f(z)=\int_{\Omega} 
K(z,w)f(w)dV(w) \quad \mbox{ for $f\in L^2(\Omega)$}, 
$$ 
where $K:\Omega\times\Omega\to{\bf C}$ is the {\it Bergman kernel} 
of the domain $\Omega$ and $dV$ is the Lebesgue measure on 
$\Omega$. 
In this paper we restrict the Bergman kernel 
on the diagonal of the domain and 
study the boundary behavior of $K(z)=K(z,z)$. 

%\vspace{.5em} 

Although there are many explicit computations  
for the Bergman kernels of specific domains 
(\cite{ber},\cite{cha},\cite{ise},\cite{dan1},\cite{grs},\cite{bol},\cite{dan3},\cite{frh1},\cite{frh2},\cite{joe1}),
it seems difficult to express the Bergman kernel 
in closed form in general.   
Therefore appropriate approximation formulas are necessary  
to know the boundary behavior of the Bergman kernel.  
From this viewpoint the following studies have great success 
in the case of strictly pseudoconvex domains. 
Assume $\Omega$ is a bounded strictly pseudoconvex domain. 
L. H\"ormander \cite{hor} shows that the limit of 
$K(z)d(z\!-\!z^0)^{n+1}$ at $z^0 \in \partial\Omega$ 
equals the determinant of 
the Levi form at $z^0$ times $n!/4\pi^n$, where $d$ is the 
Euclidean distance. 
Moreover C. Fefferman \cite{fef} and L. Boutet de Monvel and 
J. Sj\"ostrand \cite{bos} give the following asymptotic expansion 
of the Bergman kernel 
of $\Omega$ :  
\begin{equation}
K(z)=\frac{\varphi(z)}{r(z)^{n+1}}+\psi(z)\log r(z), 
\label{eqn:1.1}
\end{equation} 
where $r\in C^{\infty}(\bar{\Omega})$ is a defining function 
of $\Omega$ (i.e. $\Omega=\{r>0\}$ and  $|dr|>0$ 
on $\partial\Omega$) 
and  $\varphi$, $\psi \in C^{\infty}(\bar{\Omega})$ 
can be expanded asymptotically 
with respect to $r$. 

%\vspace{.5em} 

On the other hand,  
there are not so strong results in the weakly pseudoconvex case. 
Let us recall important studies in this case. 
Many sharp estimates of 
the size of the Bergman kernel are obtained 
(\cite{her1},\cite{ohs1},\cite{dho},\cite{cat},\cite{mcn1},   
\cite{her2},\cite{dih},\cite{her3},\cite{mcn2},\cite{cho},\cite{ohs2},\cite{goz},\cite{mcn3}). 
In particular D. Catlin \cite{cat} gives a 
complete estimate from above and below   
for domains of finite type in ${\bf C}^2$. 
Recently H. P. Boas, E. J. Straube and J. Yu \cite{bsy} have 
computed a boundary limit in the sense 
of H\"ormander for a large class of domains of finite type 
on a non-tangential cone. 
However asymptotic formulas 
are yet to be explored more extensively.  
In this paper we give an asymptotic expansion 
of the Bergman kernel for certain class of  
weakly pseudoconvex tube domains of finite type in ${\bf C}^2$. 
N. W. Gebelt \cite{geb} and F. Haslinger \cite{has2} 
have recently computed 
for the special cases, but 
the method of our expansion is different from theirs.   

%\vspace{.5em} 

Our main idea used to analyse  the Bergman kernel 
is to introduce certain real blowing-up. 
Let us briefly indicate how this blowing-up 
works for the Bergman kernel at a weakly pseudoconvex point $z^0$.  
Since the set of strictly pseudoconvex points are 
dense on the boundary of the domain of finite type,   
it is a serious problem to resolve the difficulty caused 
by strictly pseudoconvex points near $z^0$. 
This difficulty can be avoided by   
restricting the argument on a non-tangential cone 
in the domain 
(\cite{her1},\cite{dho},\cite{her2},\cite{dih},\cite{bsy}). 
We surmount the difficulty in the case of certain class of 
tube domains in the following. 
By blowing up at the weakly pseudoconvex point $z^0$,  
we introduce two new variables. 
The Bergman kernel can be developed asymptotically 
in terms of these variables  
in the sense of Sibuya \cite{sib}. 
(See also Majima \cite{maj}.)  
The expansion, regarded as a function of the first variable,  
has the form of Fefferman's expansion (\ref{eqn:1.1}), 
and hence it reflects the strict pseudoconvexity. 
The characteristic influence of the weak pseudoconvexity 
appears in the expansion with respect to 
second variable. 
Though the form of this expansion is similar to (\ref{eqn:1.1}), 
we must use $m$th root of the defining function, i.e. 
$r^{\frac{1}{m}}$, 
as the expansion variable  
when  $z^0$ is of type $2m$. 
We remark that 
a similar situation occurs in the case of another class of 
domains in \cite{geb}. 

%\vspace{.5em}

Our method of the computation is based on the studies 
\cite{fef},\cite{bos},\cite{boc},\cite{nak}. 
Our starting point is certain  integral representation in 
\cite{kor},\cite{nak}. 
After introducing the blowing-up to this representation, 
we compute the asymptotic expansion by using 
the stationary phase method. 
For the above computation, it is necessary    
to localize the Bergman kernel near a weakly pseudoconvex point.  
This localization can be obtained in a fashion similar to 
the case of some class 
of Reinhardt domains (\cite{boc},\cite{nak}). 

%\vspace{.8em}

This paper is organized as follows. 
Our main theorem is established in Section 2. 
The next three sections prepare the proof of 
the theorem. 
First an integral representation is introduced,  
which is a clue to our analysis in Section 3. 
Second the usefulness of our blowing-up 
is shown by using a simple tube domain 
$\{(z_1,z_2)\in {\bf C}^2 ; {\rm Im}z_2>[{\rm Im}z_1]^{2m}\},$ 
$m=2,3,\ldots$, in Section 4. 
This domain is considered to be a model domain for more 
general case. 
Third a localization lemma is established in Section 5, 
which is necessary to the computation in the proof of 
our theorem. 
Our main theorem is proved in Section 6. 
After an appropriate localization (\S 6.1)  
and the blowing-up at a weakly pseudoconvex point, 
an easy computation shows that 
certain two propositions are sufficient to 
prove our theorem (\S 6.2). 
In order to prove these propositions, 
we compute the asymptotic expansion of two functions 
by using the stationary phase method (\S 6.3, 6.4). 
The rest of Section 6 (\S6.5, 6.6) is devoted to 
proving two propositions. 
In Section 7 an analogous theorem about the Szeg\"o kernel 
is established. 

%\vspace{.7em}

 I would like to thank Professors Kazuo Okamoto, Takeo Ohsawa and Katsunori 
Iwasaki for generosity and several useful conversations. 
I would also like to Professor Iwasaki 
who carefully read the manuscript 
and supplied many corrections.

%%%%%%%%%%%%%%%%%%%%%%%%%%%%%%%%%%%%%%%%%%%%%%%%%%%%%%%%%%%%%%

\section{Statement of main result}

Given a function $f \in C^{\infty}({\bf R})$ 
satisfying that  
\begin{eqnarray}       
\left\{
 \begin{array}{rl}
 \!\!& \mbox{
       $f''\geq 0$ on ${\bf R}$ and $f$ has the form 
       in some neighborhood of $0$:}\\ 
 \!\!& \mbox{                                                           
       $f(x)=x^{2m}g(x)$ where $m=2,3,\ldots$, 
       $g(0)>0$ and $xg'(x) \leq 0$}.  
 \end{array}\right.  
\label{eqn:2.1}
\end{eqnarray}
Let $\omega_f \subset {\bf R}^2$ be a domain defined by  
$
\omega_f = \{(x,y) ;
 y > f(x) \}. 
$
Let $\Omega_f \subset {\bf C}^2$ 
be the tube domain over  
$\omega_f$, i.e., 
$$
\Omega_f={\bf R}^2 + i \omega_f. 
$$
Let $\pi: {\bf C}^2 \to {\bf R}^2$ be the projection defined  
by $\pi(z_1,z_2)= ({\rm Im}z_1,{\rm Im}z_2)$. 
It is easy to check that 
$\Omega_f$ is a pseudoconvex domain ; moreover 
$z^0 \in \partial\Omega_f$, with $\pi(z^0)=O$, 
is a weakly pseudoconvex point of type $2m$ (or $2m-1$)  
in the sense of Kohn or D'Angelo and 
$\partial\Omega_f \setminus \pi^{-1}(O)$ 
is strictly pseudoconvex near $z^0$. 

Now we introduce the transformation $\sigma$, 
which plays a key role on our analysis. 
Set $\Delta=\{(\tau,\varrho); 0<\tau \leq 1, \,\, \varrho>0 \}$.  
The transformation 
$\sigma : \overline{\omega_f} \to \overline{\Delta}$ is defined by 
\begin{eqnarray}       
\sigma :  
\left\{
 \begin{array}{rl}
 \!\!\!\!\! &\tau = \chi(1-\frac{f(x)}{y}), \\
 \!\!\!\!\! &\varrho= y, 
 \end{array}\right.  
\label{eqn:2.2}
\end{eqnarray}
where the function $\chi \in C^{\infty}([0,1))$ 
satisfies the conditions: 
$\chi'(u) \geq 1/2$ on $[0,1]$, and 
$\chi(u)=u$ for $0 \leq u \leq 1/3$ and 
$\chi(u)=1-(1-u)^{\frac{1}{2m}}$ 
for $1-1/3^{2m} \leq u \leq 1$.   
%Note that $\chi^{-1} \in C^{\infty}([0,1])$. 
Then $\sigma\circ\pi$ is the transformation from 
$\overline{\Omega}$ to $\overline{\Delta}$. 

The transformation $\sigma$ induces an isomorphism of 
$\omega_f \cap \{x\geq 0\}$ (or 
$\omega_f \cap \{x\leq 0\}$) on to 
$\Delta$. 
The boundary of $\omega_f$ is transfered by $\sigma$ in the following: 
$\sigma((\partial \omega_f)\setminus \{O\})= 
\{(0,\varrho); \varrho>0\}$ and 
$
\sigma^{-1}(\{(\tau,0);0\leq \tau \leq 1\})=\{O\}
$. 
This indicates that $\sigma$ is the real blowing-up 
of $\partial\omega_f$ at $O$, 
so we may say that $\sigma\circ\pi$ is the blowing-up at the  
weakly pseudoconvex point $z^0$. 
Moreover $\sigma$ patches the coordinates 
$(\tau,\varrho)$ on $\omega_f$, 
which can be considered as 
the polar coordinates around $O$. 
We call $\tau$ the angular variable 
and $\varrho$ the radial variable, respectively. 
Note that if $z$ approaches some strictly (resp. weakly) 
pseudoconvex points, 
$\tau(\pi(z))$ (resp. $\varrho(\pi(z))$) 
tends to $0$ on the coordinates $(\tau,\varrho)$.  

The following theorem asserts that the singularity  
of the Bergman kernel of $\Omega_f$ at $z^0$, 
with $\pi(z^0)=O$, can be essentially expressed in terms of 
the polar coordinates $(\tau,\varrho)$.  
%%%%%%%%%%%%%%%%%%%%
%%%%%%%%%%%%%%%%%%%%

\begin{thm}
The Bergman kernel of $\Omega_f$ has the form in some  
neighborhood of  $z^0$: 
\begin{equation}
K(z)=\frac{\Phi(\tau, \varrho^{\frac{1}{m}})}
{\varrho^{2+\frac{1}{m}}} 
+\tilde{\Phi}(\tau,\varrho^{\frac{1}{m}}) \log \varrho^{\frac{1}{m}},
\label{eqn:2.3}
\end{equation}
where 
$\Phi \in C^{\infty}((0,1]\times [0,\varepsilon))$ and 
$\tilde{\Phi} \in C^{\infty}([0,1]\times [0,\varepsilon))$   
with some $\varepsilon>0$. 

Moreover 
$\Phi$ is written in the form on the set 
$\{\tau> \alpha \varrho^{\frac{1}{2m}}\}$ with 
some $\alpha>0$: 
for every nonnegative integer $\mu_0$  
\begin{equation}
\Phi(\tau,\varrho^{\frac{1}{m}}) 
=\sum_{\mu=0}^{\mu_0} c_{\mu}(\tau) 
\varrho^{\frac{\mu}{m}} 
+ R_{\mu_0}(\tau, \varrho^{\frac{\mu}{m}})
\varrho^{\frac{\mu_0}{m}+\frac{1}{2m}}, 
\label{eqn:2.4}
\end{equation}
where 
\begin{equation}
c_{\mu}(\tau)
= \frac{\varphi_{\mu}(\tau)}
{\tau^{3+2\mu}} 
+\psi_{\mu}(\tau)\log \tau, 
\label{eqn:2.5}
\end{equation} 
for 
$\varphi_{\mu},\psi_{\mu} \in 
C^{\infty}([0,1])$,  $\varphi_0$ is positive on $[0,1]$ 
and $R_{\mu_0}$ satisfies 
$
|R_{\mu_0}(\tau,\varrho^{\frac{1}{m}})|
 \leq C_{\mu_0}[\tau-\alpha\varrho^{\frac{1}{2m}}]^{-4-2\mu_0}  
$ 
for  some positive constant $C_{\mu_0}$. 
\end{thm}
%%%%%%%%%%%%%%%%%%%%%%%%%%
%%%%%%%%%%%%%%%%%%%%%%%%%%

Let us describe the asymptotic expansion of the Bergman kernel $K$  
in more detail. 
Considering the meaning of the variables $\tau,\varrho$, 
we may say that each expansion with respect to $\tau$ or  
$\varrho^{\frac{1}{m}}$ is induced by the strict 
or weak pseudoconvexity, respectively.  
Actually the expansion (\ref{eqn:2.5}) 
has the same form as that of Fefferman 
(\ref{eqn:1.1}). 
By (\ref{eqn:2.4}),(\ref{eqn:2.5}), 
in order to see the characteristic influence 
of the weak pseudoconvexity on the singularity 
of the Bergman kernel $K$,  
it is sufficient to argue about $K$ on the region 
$$
{\cal U}_{\alpha}=\{z \in {\bf C}^2; \tau\circ\pi(z) 
> \alpha^{-1}\} 
\quad (\alpha>1). 
$$ 
This is because ${\cal U}_{\alpha}$ is the {\it widest} region 
where the coefficients $c_{\mu}(\tau)$'s are bounded. 
We call ${\cal U}_{\alpha}$ an admissible approach region of  
the Bergman kernel of $\Omega_f$ at $z^0$. 
The region ${\cal U}_{\alpha}$ seems deeply connected with 
the admissible approach regions studied 
in \cite{kra1},\cite{kra2},\cite{ala},\cite{kra3}, etc. 
We remark that on the region  ${\cal U}_{\alpha}$, 
the exchange of the expansion variable $\varrho^{\frac{1}{m}}$ 
for $r^{\frac{1}{m}}$, 
where $r$ is a defining function of $\Omega_f$ 
(e.g. $r(x,y)=y-f(x)$), 
gives no influence on the form of the expansion  
on the region ${\cal U}_{\alpha}$.

Now let us compare the asymptotic expansion (\ref{eqn:2.3}) 
on ${\cal U}_{\alpha}$ with Fefferman's expansion (\ref{eqn:1.1}).  
The essential difference between them only appears in the expansion 
variable 
(i.e. $r^{\frac{1}{m}}$ in (\ref{eqn:2.3}) and $r$ in (\ref{eqn:1.1})). 
A similar phenomenon occurs  
in subelliptic estimates for the $\bar{\partial}$-Neumann problem. 
As is well-known,  
the finite-type condition is 
equivalent to the condition that 
a subelliptic estimate holds, i.e., 
$$
|||\phi |||_{\epsilon}^2 \leq 
C(||\bar{\partial}\phi||^2+||\bar{\partial}^* \phi||^2 
+||\phi||^2) \quad\, (\epsilon>0),  
$$
(refer to \cite{koh2} for the details). 
Here,  in two dimensional case, 
this estimate holds for any $0<\epsilon\leq \frac{1}{2}$ in the 
strictly pseudoconvex case and for  $0<\epsilon\leq \frac{1}{2m}$ 
in the weakly pseudoconvex case of type $2m$, respectively.  
The difference between these two cases only appears in the value of 
$\epsilon$. 
From this viewpoint, our expansion (\ref{eqn:2.3}) 
seems to be  a natural generalization of Fefferman's expansion 
(\ref{eqn:1.1})  
in the strictly pseudoconvex case. 

%\vspace{.7em}

{\it Remarks}\,\,{\it 1}.\,\,
The idea of the blowing-up $\sigma$ is originally introduced 
in the study of the Bergman kernel of the domain 
${\cal E}_m=\{z\in {\bf C}^n; \sum_{j=1}^n |z_j|^{2m_j}<1\}$ 
$(m_j \in {\bf N}$, $m_n \neq 1)$ in \cite{joe1}. 
Since ${\cal E}_m$ has high homogeneity, 
the asymptotic expansion with respect to the radial variable does not 
appear (see also \S4). 
%\vspace{-1.1em}

{\it 2}. \,\,
If we consider the Bergman kernel on the region ${\cal U}_{\alpha}$, 
then  we can remove the condition $xg'(x)\leq 0$ in (\ref{eqn:2.1}). 
Namely even if the condition $xg'(x)\leq 0$ is not satisfied,    
we can still obtain (\ref{eqn:2.3}),(\ref{eqn:2.4}) in the theorem 
where $c_{\mu}$'s are bounded on ${\cal U}_{\alpha}$. 
But the condition $xg'(x) \leq 0$ is necessary to obtain  
the asymptotic expansion with respect to $\tau$. 
%\vspace{.5em}

{\it 3}. \,\, 
From the definition of asymptotic expansion of 
functions of several variables in \cite{sib},\cite{maj}, 
the expansion in the theorem is not complete. 
In order to get a complete asymptotic expansion, 
we must take a further blowing-up at the point 
$(\tau, \varrho)=(0,0)$.   
The real blowing-up 
$(\tau, \varrho)\mapsto (\tau, \varrho\tau^{-2m})$  
is sufficient for this purpose. 
%\vspace{.5em}

{\it 4}. \,\, 
The limit of $\varrho^{2+\frac{1}{m}}K(z)$ at $z^0$ is 
$c_0(\tau)$, so  
the boundary limit depends on the angular variable $\tau$.   
But this limit is determined uniquely ($c_0 (1)=\varphi_0(1)$)   
on a non-tangential cone in $\Omega_f$ (see \cite{bsy}).

%\vspace{.9em}

{\it Notation}.\,\,\,
In this paper we use $c$, $c_j$, or $C$ for various constants 
without further comment.

%%%%%%%%%%%%%%%%%%%%%%%%%%%%%%%%%%%%%%%%%%%%%%%%%%%%%%%%%%%%%%%%

\section{Integral representation}

In this section we give an integral representation of 
the Bergman kernel, 
which is a clue to our analysis. 
Kor\'anyi \cite{kor}, Nagel \cite{nag} and 
Haslinger \cite{has1} obtain similar representations of 
Bergman kernels or Szeg\"o kernels for certain tube domains. 

%\vspace{.5 em}
In this section we assume that  
$f\in C^{\infty}({\bf R})$ is a function 
such that $f(0)=0$ and $f''(x) \geq 0$. 
The tube domain $\Omega_f \subset {\bf C}^2$ is defined as in 
Section 2. 
Let $\Lambda, \Lambda^{\ast} \subset {\bf R}^2$ be the cones defined 
by 
\begin{eqnarray*}
\Lambda 
\!\!\!&=&\!\!\! 
\{(x,y); \mbox{$(tx,ty) \in \omega_f$ for any $t>0$} \},\\ 
\Lambda^{\ast} 
\!\!\!&=&\!\!\! 
\{(\zeta_1,\zeta_2); 
            \mbox{$x\zeta_1+y\zeta_2>0$ for any $(x,y) \in \Lambda$} \},  
\end{eqnarray*}
respectively. 
We call $\Lambda^{\ast}$ the dual cone of $\omega_f$. 
Actually $\Lambda^{\ast}$ can be  computed explicitly:  
$$
\Lambda^{\ast}=
\{(\zeta_1,\zeta_2); -R^- \zeta_2 < \zeta_1 < R^+ \zeta_2 \}, 
$$
where 
$
(R^{\pm})^{-1}=\lim_{x\to \mp\infty} f(x)|x|^{-1}>0, 
$ respectively. 
We allow that $R^{\pm}=\infty$.  
If $\lim_{|x|\to\infty} f(x)|x|^{-1-\varepsilon}>0$ 
with some $\varepsilon>0$, then $R^{\pm}=\infty$, i.e. 
$\Lambda^{\ast}=\{(\zeta_1,\zeta_2);\zeta_2>0\}$. 

The Bergman kernel of $\Omega$ is expressed in the following. 
Set $(x,y)=({\rm Im} z_1,{\rm Im} z_2)$. 
\begin{equation}
K(z)=\frac{1}{(4\pi)^2} \int\!\!\!\int_{\Lambda^{\ast}} 
e^{-x\zeta_1-y\zeta_2} 
\frac{\zeta_2}{D(\zeta_1,\zeta_2)} d\zeta_1d\zeta_2, 
\label{eqn:3.3} 
\end{equation}
where 
\begin{equation}
D(\zeta_1,\zeta_2)
=\int_{-\infty}^{\infty} e^{-\xi\zeta_1-f(\xi)\zeta_2} d\xi.
\label{eqn:3.4}
\end{equation}
The above representation  
can be obtained by a slight generalization 
of the argument of Kor\'anyi \cite{kor}, 
so we omit the proof.

%%%%%%%%%%%%%%%%%%%%%%%%%%%%%%%%%%%%%%%%%%%%%%%%%%%%%%%%%%%%%%%

\section{Analysis on a model domain}

Let $\omega_0 \subset {\bf R}^2$ be a domain defined by 
$\omega_0=\{(x,y); y>g x^{2m}\}$, 
where $m =2,3,\ldots$ and $g>0$. 
Set $\Omega_0={\bf R}^2+i \omega_0$. 
F. Haslinger \cite{has2} computes the asymptotic expansion 
of the Bergman kernel of $\Omega_0$ (not only on the diagonal  
but also off the diagonal).  
In his result Fefferman's expansion only appears.  

In this paper,    
we consider       
$\Omega_0$ as a model domain 
for the study of singularity of the Bergman kernel 
for more general domains. 
The following proposition shows the reason why 
we take $\Omega_0$ as a model domain. 
Set $(x,y)=({\rm Im}z_1, {\rm Im}z_2)$. 
%%%%%%%%%%%%%%%%%%%%%%%%%%%%%%%%%%%%%%%%%%%%%%%%%%%%%%%%%%
\begin{proposition}
The Bergman kernel $K$ of $\Omega_0$ has the form: 
\begin{equation}
K(z)=\frac{\Phi(\tau)}{\varrho^{2+\frac{1}{m}}}, 
\label{eqn:4.1}
\end{equation} 
where $\tau= \chi(g x^{2m} y^{-1})$, 
$\varrho=y$ $(see\,\,\, (\ref{eqn:2.2}))$ and 
$$
\Phi(\tau)=\frac{\varphi(\tau)}{\tau^3} + 
\psi(\tau)\log \tau, 
$$
with 
$\varphi, \psi \in C^{\infty}([0,1])$ and $\varphi$ is positive on 
$[0,1]$. 
\end{proposition}
%%%%%%%%%%%%%%%%%%%%%%%%%%%%%%%%%%%%%%%%%%%%%%%%%%%%%%%%%%

{\it Proof}. \,\,\, 
Normalizing the integral representation (\ref{eqn:3.3}) and 
introducing the variables 
$t= g^{\frac{1}{2m}} x y^{-\frac{1}{2m}}$, 
$\varrho=y^{\frac{1}{2m}}$, 
we have (\ref{eqn:4.1}) 
where 
\begin{eqnarray}
&&\Phi(\tau)= \frac{2m}{(4\pi)^2} 
g^{\frac{1}{m}} \int_0^{\infty} 
e^{-s^{2m}} L(ts) s^{4m+1} ds, 
\label{eqn:4.3}\\ 
&&L(u)=\int_{-\infty}^{\infty} 
e^{uv}\frac{1}{\phi(v)}dv,\nonumber\\
&&\phi(v)=\int_{-\infty}^{\infty} 
e^{-w^{2m}+vw}dv. \nonumber
\end{eqnarray}
It turns out 
from (\ref{eqn:4.3}) and the definition of $\tau$ 
that  $\Phi \in C^{\infty}((0,1])$. 

Now let $\hat{\Phi}$ be defined by 
\begin{equation}
\hat{\Phi}(\tau)
=\frac{2m}{(4\pi)^2}
\int_1^{\infty} e^{-s^{2m}} L(ts)s^{4m+1} ds. 
\label{eqn:4.5}
\end{equation}
If we admit Lemma 6.2 in Subsection 6.4  below,   
we have 
\begin{equation} 
L(u)=u^{2m-2} e^{u^{2m}} \tilde{L}(u), 
\label{eqn:4.6}
\end{equation}
where $\tilde{L}(u)\sim \sum_{j=0}^{\infty}c_j u^{-2mj}$ 
as $u\to\infty$. 
Substituting (\ref{eqn:4.6}) into (\ref{eqn:4.5}), 
we have 
\begin{eqnarray}
\hat{\Phi}(\tau)
\!\!\!&=&\!\!\! 
\frac{2m}{(4\pi)^2} t^{2m-2} \int_1^{\infty} 
e^{-[1-t^{2m}]s^{2m}} \tilde{L}(ts)s^{6m-1} ds
\nonumber\\
\!\!\!&=&\!\!\! 
\frac{1}{(4\pi)^2} \int_1^{\infty} 
e^{-\chi^{-1}(\tau)\sigma} \hat{L}(\tau,\sigma)\sigma^2 d\sigma.
\nonumber
\end{eqnarray}
Since $\hat{L}(\tau,\sigma)\sim 
\sum_{j=0}^{\infty}c_j(\tau) \sigma^{-j}$ 
as $\sigma \to \infty$ for $c_j \in C^{\infty}([0,1])$, 
we have 
$$
\hat{\Phi}(\tau)=\frac{\hat{\varphi}(\tau)}{\tau^3} + 
\hat{\psi}(\tau)\log \tau, 
$$ 
with 
$\hat{\varphi}, \hat{\psi} \in C^{\infty}([0,1])$ and 
$\hat{\varphi}$ is positive on $[0,1]$. 

Finally since the difference between 
$\Phi$ and $\hat{\Phi}$ is smooth on $[0,1]$,  
we can obtain Proposition 4.1. 
\qed

%%%%%%%%%%%%%%%%%%%%%%%%%%%%%%%%%%%%%%%%%%%%%%%%%%%%%%%%%%%

\section{Localization lemma}

In this section we prepare a lemma, 
which  is necessary for the proof of Theorem 2.1.
This lemma shows that the singularity of the 
Bergman kernel for certain class of domains is 
determined by the local information about 
the boundary. 
The method of the proof is similar     
to the case of some class of Reinhardt domains 
(\cite{boc},\cite{nak}). 
Throughout this section, 
$j$ stands for $1$ or $2$. 
%\vspace{.4 em}

Let $f_1$, $f_2 \in C^{\infty}({\bf R})$ be functions such that 
$f_j(0)=f'_j(0)=0$, $f_j''\geq 0$ on ${\bf R}$ 
and $f_1(x)=f_2(x)$ on $|x|<\delta$.  
Let $\omega_j \subset {\bf R}^2$ be a domain  defined by 
$\omega_j =\{(x,y): y>f_j(x)\}$.  
Set $\Omega_j={\bf R}^2+i \omega_j \subset {\bf C}^2$.

%
%%%%%%%%%%%%%%%%%%%%%%%%%%%%%%%%%%%%%%
\begin{lemma} 
Let $K_j$ be the Bergman kernels 
of $\Omega_j$ for $j=1,2$, respectively. 
Then we have  
$$
K_1(z)- K_2(z) \in C^{\omega}(U),  
$$
where $U$ is some neighborhood of $z^0$. 

\end{lemma}
%%%%%%%%%%%%%%%%%%%%%%%%%%%%%%%%%%%%%%
%

{\it Proof}. \,\,\, 
Let $\Lambda_j^{\ast}$  be the dual cone  
of $\omega_j$, i.e. $\Lambda_j^{\ast} 
=\{(\zeta_1,\zeta_2); -R_j^- \zeta_2 < \zeta_1 <  R_j^+ \zeta_2 \}$, 
where $(R_j^{\pm})^{-1}=\lim_{x\to\mp} f(x)|x|^{-1}$, 
respectively (see \S3). 
Let $K_j[\Delta](x,y)$ be defined by 
$$
K_j[\Delta](x,y)
=\frac{1}{(4\pi)^2}\int\!\!\!\int_{\Delta} e^{-y \zeta_2 -x \zeta_1} 
\frac{\zeta_2}{D_j(\zeta_1,\zeta_2)} d\zeta_1 d\zeta_2,
$$  
where $\Delta \subset {\bf R}^2$ and 
$D_j(\zeta_1,\zeta_2)
=\int_{-\infty}^{\infty} e^{-\zeta_2 f_j(\xi) - \zeta_1\xi} d\xi$. 
Set $\Lambda_{\varepsilon} 
=\{ (\zeta_1,\zeta_2); |\zeta_1|< \varepsilon\zeta_2 \}$, where   
$\varepsilon>0$ is small. 
Now the following claims (i), (ii) imply Lemma 5.1. 
Set $O=(0,0)$.  
\begin{eqnarray*}
&{\rm (i)}& K_j[\Lambda_j^{\ast}] \equiv 
K_j[\Lambda_{\varepsilon}]  \quad \mbox{modulo $C^{\omega}(\{O\})$ 
\,\, for any $\varepsilon>0$}, \\
&{\rm (ii)}& K_1[\Lambda_{\varepsilon_0}] \equiv 
K_2[\Lambda_{\varepsilon_0}] \quad 
\mbox{modulo $C^{\omega}(\{O\})$ for some $\varepsilon_0>0$}. 
\end{eqnarray*}
In fact if we substitute $(x,y)\!=\!({\rm Im}z_1,{\rm Im}z_2)$, then 
$K_1\!=\!K_1[\Lambda_1^{\ast}] \!\equiv\! K_1[\Lambda_{\varepsilon_0}]
\!\equiv\! K_2[\Lambda_{\varepsilon_0}]
\!\equiv\! K_2[\Lambda_2^{\ast}] 
\!=\! K_2$ modulo $C^{\omega}(\{z^0\})$. 

%\vspace{.5em}
Let us show the above claims. 
(i)\,\,\, 
Set 
$\Lambda_{\varepsilon}^{\pm} 
=\{(\zeta_1,\zeta_2); 0 < \varepsilon\zeta_2 < \pm\zeta_1 
<  R_j^{\pm} \zeta_2 \}$, 
respectively.
Since $K_j[\Lambda_j^{\ast}]-K_j[\Lambda_{\varepsilon}]
=K_j[\Lambda_{\varepsilon}^+]+K_j[\Lambda_{\varepsilon}^-]$, 
it is sufficient to show 
$K_j[\Lambda_{\varepsilon}^{\pm}] \in C^{\omega}(\{O\})$.  
We only consider the case of $K_j[\Lambda_{\varepsilon}^{+}]$. 
Changing the integral variables, 
we have 
$$
K_j[\Lambda_{\varepsilon}^+](x,y)
=\frac{1}{(4\pi)^2} 
\int_0^{\infty}\!\!\int_{\varepsilon}^{R_j^+}
H_j(\zeta,\eta; x,y) d\zeta d\eta, 
$$
where 
\begin{eqnarray*}
&&H_j(\zeta,\eta; x,y)\,(=H_j) = 
e^{-y\eta+x\eta\zeta}
\frac{\eta^2}{E_j(\zeta,\eta)}, \\
&&E_j(\zeta,\eta)=
\int_{-\infty}^{\infty}
e^{-\eta[f_j(\xi)-\zeta\xi]}d\xi.
\end{eqnarray*}
It is an important remark that 
$K_j[\Lambda_{\varepsilon}^+]$ is real analytic 
on the region where $H_j$ is integrable on 
$\{(\zeta,\eta); \zeta> \varepsilon, \eta>0\}$. 

If we take $\delta_1>0$ such that 
$|f_j(\xi)\xi^{-1}|<\frac{1}{2}\varepsilon$ if $|\xi|<\delta_1$,  
then we have 
\begin{eqnarray}
E_j(\zeta,\eta)
&\geq& 
2\int_0^{\delta_1} e^{\eta\xi[\zeta-f(\xi)\xi^{-1}]}d\xi \nonumber\\ 
&\geq& 
2\int_0^{\delta_1} e^{\frac{1}{2}\varepsilon\eta\xi}d\xi 
\geq
\frac{4}{\varepsilon\eta}[e^{\frac{1}{2}\delta_1\varepsilon\eta}-1]
\label{eqn:5.4} 
\end{eqnarray}
By (\ref{eqn:5.4}), 
we have $H_j(\zeta,\eta;x,y)\leq C \eta 
e^{-[y-x\zeta+\frac{1}{2}\delta_1\varepsilon]\eta}$ for $\eta\geq 1$. 
This inequality implies that 
if $x<0$ and $y>-\frac{1}{2}\delta_1\varepsilon$, 
then $H_j$ is integrable 
on $\{(\zeta,\eta); \zeta> \varepsilon, \eta>0\}$. 
Thus $K_j[\Lambda_{\varepsilon}^+]$ is 
real analytic on $\omega_j \cup 
\{(x,y); x<0, y>-\frac{1}{2}\delta_1\varepsilon\} 
=:\omega_j^+$. 
By regarding $x,y$ as two complex variables, 
$K_j[\Lambda_{\varepsilon}^+](x,y)$  
is holomorphic on $\omega_j^+ +i{\bf R}^2$, 
so the shape of $\omega_j^+$ implies that 
$K_j[\Lambda_{\varepsilon}^+]$ can be extended holomorphically 
to a region containing some neighborhood of $\{O\}+i {\bf R}^2$. 
Consequently 
we have $K_j[\Lambda_{\varepsilon}^+]\in C^{\omega}(\{O\})$.

%%%%%%%%%%%%%%%

%\vspace{.5 em}

(ii)\,\,\, 
Changing the integral variables, we have 
$$
K_1[\Lambda_{\varepsilon}](x,y)- 
K_2[\Lambda_{\varepsilon}](x,y) = 
\frac{1}{(4\pi)^2}\int_{0}^{\infty}\!\!\!
\int_{-\varepsilon}^{\varepsilon}
(H_1-H_2) d\zeta d\eta.
$$
We remark that 
$K_1[\Lambda_{\varepsilon}]- K_2[\Lambda_{\varepsilon}]$ is 
real analytic on the region where $H_1-H_2$ is integrable on 
$\{(\zeta,\eta); |\zeta|<\varepsilon, \eta>0\}$.   
To find a positive number $\varepsilon_0$ satisfying (ii), 
let us consider the integrability of 
\begin{equation}
|H_1-H_2|
=\eta^2 e^{-y\eta+x\eta\zeta}
\frac{|E_2(\zeta,\eta)-E_1(\zeta,\eta)|}
{|E_1(\zeta,\eta) \cdot E_2(\zeta,\eta)|}. 
\label{eqn:5.41}
\end{equation}

First we give an estimate of  
$|E_2(\zeta,\eta)-E_1(\zeta,\eta)|$. 
Let $\varepsilon_1 >0$ be defined by 
$|f_j(\xi)\xi^{-1}|\geq \varepsilon_1$ if $|\xi| \geq \delta$. 
If $|\zeta|\leq \frac{1}{2}\varepsilon_1$, then 
\begin{equation}
\int_{|\xi|\geq \delta} 
e^{-\eta\xi[f_j(\xi)\xi^{-1}-\zeta]} d\xi 
\leq  2 \int_{\delta}^{\infty}  
e^{-\frac{1}{2}\varepsilon_1\eta\xi} d\xi \,\,\,
\leq \frac{4}{\varepsilon_1\eta}
e^{-\frac{1}{2}\delta\varepsilon_1\eta}. 
\label{eqn:5.5} 
\end{equation}
Therefore (\ref{eqn:5.5}) implies 
\begin{equation}
|E_2(\zeta,\eta)-E_1(\zeta,\eta)| 
\leq 
\sum_{j=1}^2 \int_{|\xi|\geq \delta} 
e^{-\eta[f_j(\xi)-\zeta\xi]} d\xi
\leq 
\frac{8}{\varepsilon_1 \eta}
e^{-\frac{1}{2}\delta\varepsilon_1\eta}.
\label{eqn:5.6} 
\end{equation}

Second we give an estimate of 
$E_j(\zeta,\eta)$. 
By Taylor's formula, 
we can choose $\varepsilon_2>0$ satisfying the following. 
If $|\zeta|<\varepsilon_2$, 
then there is a function 
$\alpha_j(\zeta)=\alpha_j$ ($\alpha_j(0)=0$)  
such that $f_j'(\alpha_j)=\zeta$ and moreover 
there is a bounded function $R_j(\zeta,\xi)$ on 
$[-\varepsilon_2,\varepsilon_2] \times [-\xi_0,\xi_0]$, 
with some $\xi_0>0$, such that 
$F_j(\xi)$ $(:=f_j(\xi)-\zeta\xi)$ 
$=F_j(\alpha_j)+R_j(\zeta, \xi-\alpha_j)(\xi-\alpha_j)^2$  
with $F_j(\alpha_j)\leq 0$. 
Then if $|\zeta|<\varepsilon_2$,  
we have 
\begin{eqnarray}
E_j(\zeta,\eta)
\!\!\!&=&\!\!\!
\int_{-\infty}^{\infty} e^{-\eta F_j(\xi)} d\xi 
\geq e^{-\eta F_j(\alpha_j)} 
\int_{-\xi_0}^{\xi_0} e^{-\eta R_j(\zeta,\xi)\xi^2} d\xi 
\nonumber\\ 
\!\!\!&\geq&\!\!\! 
\frac{C}{\sqrt{\eta}} e^{-\eta F_j(\alpha_j)} 
\geq \frac{C}{\sqrt{\eta}}. 
\label{eqn:5.8} 
\end{eqnarray}

Now we set $\varepsilon_0=\min\{\varepsilon_1,\varepsilon_2\}$. 
Then by putting 
(\ref{eqn:5.41}),(\ref{eqn:5.6}),(\ref{eqn:5.8}) 
together, we have 
$$
|H_1-H_2|
\leq C \eta^2 e^{-\eta[y-x\zeta+\frac{1}{2}\delta\varepsilon_1]} 
\quad \,\,\,\mbox{for $|\zeta|<\varepsilon_0, \eta>0$}.
$$
This inequality implies that 
if $\varepsilon_0 |x|-y-\frac{1}{2}\delta\varepsilon_1>0$, 
then $|H_1-H_2|$ is integrable on 
$\{(\zeta,\eta);|\zeta|<\varepsilon_0, \eta>0\}$. 
Hence $K_1[\Lambda_{\varepsilon_0}]-K_2[\Lambda_{\varepsilon_0}]$ 
is real analytic on the region 
$\{(x,y); \varepsilon_0 |x|-y-\frac{1}{2}\delta\varepsilon_1 >0\}$, 
which contains $\{O\}$.  

This completes the proof of Lemma 5.1. 
\qed

%%%%%%%%%%%%%%%%%%%%%%%%%%%%%%%%%%%%%%%%%%

%%%%%%%%%%%%%%%%%%%%%%%%%%%%%%%%%%%%%%%%%%  

\section{Proof of Theorem 2.1}

In this section we give a proof of Theroem 2.1. 
The definitions of $f$, $\omega_f$ and $\Omega_f$ are given as in 
Section 2.

\subsection{Localization}

From the previous section, it turns out that 
the singularity of the Bergman kernel of $\Omega_f$ at $z^0$ 
is determined by the local information about $\partial\Omega_f$ 
near $z^0$.  
Thus we construct an appropriate domain 
whose boudary coincides $\partial\Omega_f$ near $z^0$ 
for the computation below. 

%\vspace{.5em}
We can easily construct 
a function $\tilde{g} \in C^{\infty}({\bf R})$ 
such that 
\begin{eqnarray}       
&\tilde{g}(x)=\left\{
 \begin{array}{rl}
  g(x) & \quad \mbox{for $|x|\leq \delta$} \\
  \frac{9}{10} g(0) & \quad \mbox{for $|x|\geq 1$}
 \end{array}\right.  \quad \mbox{and}
\label{eqn:6.1} \\
&0\leq -x\tilde{g}'(x),\, |x^2\tilde{g}''(x)| < \frac{1}{5}g(0) \quad 
\mbox{for $x\in {\bf R}$}, 
\label{eqn:6.2} 
\end{eqnarray}
for some small positive constant $\delta < 1$. 
Note that $\frac{9}{10}g(0)\leq \tilde{g}(x)\leq g(0)$. 
Set $\tilde{f}(x)=x^{2m}\tilde{g}(x)$ 
and $\omega_{\tilde{f}}=
\{
(x,y)\in {\bf R}^2; y>\tilde{f}(x)
\}$. 
Let $\Omega_{\tilde{f}} \subset {\bf C}^2$ 
be the tube domain over $\omega_{\tilde{f}}$, i.e. 
$\Omega_{\tilde{f}}={\bf R}^2 + i \omega_{\tilde{f}}$. 
Here we remark that the boundary of $\Omega_{\tilde{f}}$ 
is strictly pseudoconvex off the set 
$\{(z_1,z_2);{\rm Im}z_1\!=\!{\rm Im} z_2\!=\!0\}$. 
In fact we can easily check that 
$\tilde{f}''(x)>0$ if $x\neq0$ by 
(\ref{eqn:6.1}),(\ref{eqn:6.2}). 
 
Let $\tilde{K}$ be the Bergman kernel of $\Omega_{\tilde{f}}$. 
In order to obtain Theorem 2.1, 
it is sufficient to consider the 
singularity of the Bergman kernel $\tilde{K}$ near $z^0$ 
by Lemma 5.1.

%%%%%%%%%%%%%%%%%%%%%%%%%%%%%%%%%%%%%%%%%%%%%%%%%%%%%%%%%%
\subsection{Two propositions and the proof of Theorem 2.1}

A clue to our analysis of the Bergman kernel 
is the integral representation in Section 3. 
Normalizing this representation, 
the Bergman kernel $\tilde{K}$ of $\Omega_{\tilde{f}}$ 
can be expressed in the following. 
$$
\tilde{K}(z)=
\frac{2m}{(4\pi)^2}g(0)^{\frac{1}{m}} 
\int_0^{\infty} e^{-y u^{2m}} P(x,u) u^{4m+1} du,
$$
with 
\begin{eqnarray*}
P(x,u) \!\!\!&=&\!\!\! \int_{-\infty}^{\infty}
      e^{g(0)^{\frac{1}{2m}}xuv} \frac{1}{\phi(v,u^{-1})} dv,\\
\phi(v,X)
       \!\!\!&=&\!\!\! \int_{-\infty}^{\infty}
e^{-\hat{g}(Xw)w^{2m} +vw}dw, 
\end{eqnarray*}
where $\hat{g}(x)=\tilde{g}(x)/ g(0)$. 
In order to prove Theorem 2.1, 
it is sufficient to  consider the following function $\bar{K}$  
instead of $\tilde{K}$. 
\begin{equation}
\bar{K}(z)=
\frac{2m}{(4\pi)^2} g(0)^{\frac{1}{m}} 
\int_1^{\infty} e^{-y u^{2m}} P(x,u) u^{4m+1} du.    
\label{eqn:6.5}
\end{equation}
In fact the difference between $\tilde{K}$ and $\bar{K}$ 
is smooth. 
%\vspace{.4em}

By introducing the variables 
$t_0=g(0)^{\frac{1}{2m}}x y^{\frac{-1}{2m}}$, 
$\xi=y^{\frac{1}{2m}}$ to the integral representation 
(\ref{eqn:6.5}), 
we have 
\begin{equation}
\bar{K}(z)=
\frac{2m}{(4\pi)^2}\xi^{-4m-2}g(0)^{\frac{1}{m}}
\int_{\xi}^{\infty}
e^{-s^{2m}}L(t_0,\xi;s)s^{4m+1}ds,
\label{eqn:6.051}
\end{equation}
where  
\begin{equation}
L(t_0,\xi;s)=\int_{-\infty}^{\infty}e^{t_0 sv}
\frac{1}{\phi(v, \xi s^{-1})} dv.
\label{eqn:6.052}
\end{equation} 

We divide the integral in (\ref{eqn:6.051}) into two parts: 
\begin{equation}
\bar{K}(z)=\frac{2m}{(4\pi)^2}g(0)^{\frac{1}{m}}\xi^{-4m-2} 
\{K^{\langle 1 \rangle}(z)+K^{\langle 2 \rangle}(z)\},
\label{eqn:6.053}
\end{equation}
where
\begin{eqnarray}
K^{\langle 1 \rangle}(z)
\!\!\!&=&\!\!\! 
\int_{1}^{\infty}
e^{-s^{2m}}L(t_0,\xi;s)s^{4m+1}ds,
\label{eqn:6.6}\\
K^{\langle 2 \rangle}(z)
\!\!\!&=&\!\!\!
\int_{\xi}^{1}
e^{-s^{2m}}L(t_0,\xi;s)s^{4m+1}ds.
\label{eqn:6.7} 
\end{eqnarray}

%%%%%%%%%%%%%%%%%%%%%%%%%
Since the function 
$[\phi(v,X)]^{-1}$ is smooth function of $X$ 
on $[0,1]$, 
for any positive integer $\mu_0$
\begin{equation}
\frac{1}{\phi(v,X)}
=\sum_{\mu=0}^{\mu_0} 
a_{\mu}(v)X^{\mu} + r_{\mu_0}(v,X)X^{\mu_0+1}, 
\label{eqn:6.071}
\end{equation}
where 
\begin{eqnarray}
a_{\mu}(v) 
\!\!\!& = &\!\!\! 
\left. \frac{1}{\mu!}
\frac{\partial^{\mu}}{\partial X^{\mu}} 
\frac{1}{\phi(v,X)}\right|_{X=0}, 
\label{eqn:6.8}\\
r_{\mu_0}(v,X) 
\!\!\!& = &\!\!\! 
\left. 
\frac{1}{\mu_0!} 
\int_{0}^{1}
(1-p)^{\mu_0} 
\frac{\partial^{\mu_0+1}}{\partial Y^{\mu_0+1}}
\frac{1}{\phi(v,Y)}\right|_{Y=Xp} dp.
\nonumber 
\end{eqnarray}
Substituting (\ref{eqn:6.071}) into (\ref{eqn:6.052}), 
we have 
\begin{equation}
L(t_0,\xi;s)
=\sum_{\mu=0}^{\mu_0}L_{\mu}(t_0 s)\xi^{\mu} s^{-\mu} 
+ \tilde{L}_{\mu_0}(t_0,\xi;s)\xi^{\mu_0+1}s^{-\mu_0-1},
\label{eqn:6.091}
\end{equation}
where 
\begin{eqnarray*}
L_{\mu}(u) 
\!\!\!&=&\!\!\!
\int_{-\infty}^{\infty}
e^{uv}a_{\mu}(v) dv,\\
\tilde{L}_{\mu_0}(t_0,\xi;s)
\!\!\!&=&\!\!\!
\int_{-\infty}^{\infty}
e^{t_0 sv}r_{\mu_0}(v,\xi s^{-1})dv. 
\end{eqnarray*}

Moreover substituting (\ref{eqn:6.091}) into 
(\ref{eqn:6.6}),(\ref{eqn:6.7}),  
we have 
$$
K^{\langle j \rangle}(z)=
\sum_{\mu=0}^{\mu_0} K_{\mu}^{\langle j \rangle}(\tau,\xi) 
\xi^{\mu}
+ \tilde{K}_{\mu_0}^{\langle j \rangle}(\tau,\xi) \xi^{\mu_0+1}, 
$$
for $j = 1,2$ where 
\begin{eqnarray*}
K_{\mu}^{\langle 1 \rangle}(\tau,\xi)
\!\!\!&=&\!\!\! 
\int_{1}^{\infty}
e^{-s^{2m}}L_{\mu}(t_0 s)s^{4m+1-\mu}ds,\\
\tilde{K}_{\mu_0}^{\langle 1 \rangle}(\tau,\xi)
\!\!\!&=&\!\!\! 
\int_{1}^{\infty}
e^{-s^{2m}} \tilde{L}_{\mu_0}(t_0,\xi;s)s^{4m-\mu_0}ds,\\
K_{\mu}^{\langle 2 \rangle}(\tau,\xi)
\!\!\!&=&\!\!\! 
\int_{\xi}^{1}
e^{-s^{2m}}L_{\mu}(t_0 s)s^{4m+1-\mu}ds,\\
\tilde{K}_{\mu_0}^{\langle 2 \rangle}(\tau,\xi)
\!\!\!&=&\!\!\!  
\int_{\xi}^{1}
e^{-s^{2m}} \tilde{L}_{\mu_0}(t_0,\xi;s)s^{4m-\mu_0}ds. 
\end{eqnarray*}
 
The following two propositions are concerned 
with the singularities of 
the above functions. 
Their proofs are given in Subsections 6.5, 6.6.

%%%%%%%%%%%%%%%%%%%%%%%%%%%%%%%%%%%%%%%%%%%%%%%%%%%%%%%%
\begin{proposition}

$(i)$\,\,For any nonnegative integer $k_0$, 
$K_{\mu}^{\langle 1 \rangle}$ is expressed in the form: 
$$
K_{\mu}^{\langle 1 \rangle}(\tau,\xi)
=
\sum_{k=0}^{k_0} c_{\mu,k}(\tau)\xi^k 
+ \tilde{K}_{\mu,k_0}^{\langle 1 \rangle}(\tau,\xi)\xi^{k_0+1}, 
$$
where 
$$
c_{\mu,k}(\tau)
= \frac{\varphi_{\mu,k}(\tau)}
{\tau^{3+\mu+k}} 
+\psi_{\mu,k}(\tau)\log \tau, 
$$
for $\varphi_{\mu,k},\psi_{\mu,k} \in C^{\infty}([0,1])$ 
and $\tilde{K}_{\mu,k_0}^{\langle 1 \rangle}$ satisfies 
$|\tilde{K}_{\mu,k_0}^{\langle 1 \rangle}(\tau,\xi)|<
C_{\mu,k_0}[\tau-\alpha\xi]^{-4-\mu-k_0}$  
for some positive constants $C_{\mu,k_0}$ and $\alpha$. 
 
$(ii)$\,\,$\tilde{K}_{\mu_0}^{\langle 1 \rangle}$ satisfies 
$
|\tilde{K}_{\mu_0}^{\langle 1 \rangle}(\tau,\xi)| < 
C_{\mu_0}[\tau-\alpha\xi]^{-4-\mu_0}
$
for some positive constants $C_{\mu_0}$ and $\alpha$. 
\end{proposition}
%%%%%%%%%%%%%%%%%%%%%%%%%%%%%%%%%%%%%%%%%%%%%%%%%%%%%%%%%%%%%%%

%%%%%%%%%%%%%%%%%%%%%%%%%%%%%%%%%%%%%%%%%%%%%%%%%%%%%%%%%%%%%%%
\begin{proposition}

$(i)$\,\,$(a)$\,\,For $0 \leq \mu \leq 4m+1$,
$
K_{\mu}^{\langle 2 \rangle} \in 
C^{\infty}([0,1] \times [0,\varepsilon)). 
$

$(b)$\,\, 
For $\mu \geq 4m+2$, 
$K_{\mu}^{\langle 2 \rangle}$ 
can be expressed in the form: 
$$
K_{\mu}^{\langle 2 \rangle}(\tau,\xi) \xi^{-4m-2+\mu} 
=H_{\mu}(\tau, \xi) \xi^{-4m-2+\mu} \log \xi 
+\tilde{H}_{\mu}(\tau,\xi), 
$$
where $H_{\mu},\tilde{H}_{\mu} \in C^{\infty}
([0,1] \times [0,\varepsilon))$. 

$(ii)$\,\,
For any positive integer $r$, 
there is a positive integer $\mu_0$ such that 
$$
\tilde{K}_{\mu_0}^{\langle 2 \rangle}(\tau,\xi) \xi^{-4m-1+\mu_0} 
\in 
C^{r}([0,1] \times [0,\varepsilon)).
$$
\end{proposition}

%%%%%%%%%%%%%%%%%%%%%%%%%%%%%%%%%%%%%%%%%%%%%%%%%%%%%%%%%%%%%
First by  Proposition 6.1, 
$K^{\langle 1 \rangle}$ 
can be expressed in the form: 
\begin{equation}
K^{\langle 1 \rangle}(z)
= 
\sum_{\mu=0}^{\mu_0}
c_{\mu}(\tau)\xi^{\mu} 
+R_{\mu_0}(\tau,\xi) \xi^{\mu_0+1},
\label{eqn:6.200}
\end{equation} 
where $c_{\mu}$'s are expressed as in (\ref{eqn:2.5}) 
in Theorem 2.1 and 
$R_{\mu_0}$ satisfies $|R_{\mu_0}|<C_{\mu_0}
[\tau-\alpha\xi]^{-4-\mu_0}$.

Next by Proposition 6.2, 
$K^{\langle 2 \rangle}(z)$ 
can be expressed in the form: 
for any positive integer $r$, 
\begin{equation}
\xi^{-4m-2} 
K^{\langle 2 \rangle}(z)
=H(\tau,\xi) \log \xi + \tilde{H}(\tau,\xi), 
\label{eqn:6.210}
\end{equation}
where $H\in C^{\infty}([0,1] \times [0,\varepsilon))$ 
and $\tilde{H} \in C^{r}([0,1] \times [0,\varepsilon))$.

%\vspace{.5 em}

Hence 
putting 
(\ref{eqn:6.053}),(\ref{eqn:6.200}),(\ref{eqn:6.210}) together, 
we can obtain Theorem 2.1. 
Note that $K(z)$ is an even function of $\xi$. 

\qed

%%%%%%%%%%%%%%%%%%%%%%%%%%%%%%%%%%%%%%%%%%%%%%%%%%%%%%%%%%%%%%
\subsection{Asymptotic expansion of $a_{\mu}$}

By a direct computation in (\ref{eqn:6.8}),  
$a_{\mu}(v)$ can be expressed in the following form: 
\begin{equation}
a_{\mu}(v)
=
\sum_{|\alpha|=\mu}
C_{\alpha} \frac{
\phi^{[\alpha_1]}(v)\cdots \phi^{[\alpha_{\mu}]}(v)
}{\phi(v)^{\mu+1}}, 
\label{eqn:6.16}
\end{equation}
where $\phi^{[k]}(v)=
\left.\frac{\partial^k}{\partial X^k}\phi(v,X)\right|_{X=0}$, 
$C_{\alpha}$'s are constants depending on 
$\alpha=(\alpha_1,\ldots,\alpha_{\mu})\in {\bf Z}_{\geq 0}^{\mu}$ 
and $|\alpha|=\sum_{j=1}^{\mu}\alpha_j$. 
Since 
$$
\frac{\partial^k}{\partial X^k}\phi(v,X)
=\int_{-\infty}^{\infty}
\left\{
w^k\sum_{j=1}^k c_{kj}(Xw)w^{2mj}
\right\}   
e^{-\hat{g}(Xw)w^{2m} +vw}dw,
$$
for $k \geq 1$ where 
$c_{kj} \in C^{\infty}({\bf R})$,  
we have 
\begin{equation}
\phi^{[k]}(v)=
\left.\frac{\partial^k}{\partial X^k}\phi(v,X)\right|_{X=0}
=\sum_{j=1}^k c_{kj}(0)\phi_{2mj+k}(v), 
\label{eqn:6.161}
\end{equation}
for $k \geq 1$ where 
$$\phi_{l}(v)=\int_{-\infty}^{\infty}
w^{l} e^{-w^{2m} +vw}dw.
$$
Here the following lemma 
is concerned with the asymptotic expansion of 
$\phi_{l}$ at infinity. 
%%%%%%%%%%%%%%%%%%%%%%%%%%%%%%%%%%%%%%%%%%%%%%%%%
\begin{lemma}
Set  $a=[(2m)^{\frac{-1}{2m-1}}-(2m)^{\frac{-2m}{2m-1}}]>0$. 
Then we have 
$$
\phi_l(v) \sim v^{\frac{1-m+l}{2m-1}}\cdot
\exp\{a v^{\frac{2m}{2m-1}}\}\cdot 
\sum_{j=0}^{\infty}c_j v^{-\frac{2m}{2m-1}j}
\,\,\,\,\,\,\,{\rm as}\,\,v \to +\infty\,\,\,\,
{\rm for}\,\,l \geq 0.
$$
%$$
%\frac{\varphi_{l}(v)}{\varphi(v)} \sim v^{\frac{l}{2m-1}}\cdot 
%\sum_{j=0}^{\infty}c_j v^{-\frac{2m}{2m-1}j}
%\,\,\,\,\,\,\,{\rm as}\,\,v\to+\infty.
%$$
\end{lemma}
%%%%%%%%%%%%%%%%%%%%%%%%%%%%%%%%%%%%%%%%%%%%%%%%%
The proof of the above lemma will be given soon later. 

Lemma 6.1 and (\ref{eqn:6.161}) imply 
$$
\frac{\phi^{[k]}(v)}
{\phi(v)} 
\sim 
v^{\frac{(2m+1)k}{2m-1}} \cdot
\sum_{j=0}^{\infty}c_j v^{-\frac{2m}{2m-1}j}
\,\,\,\,\,\,\,{\rm as}\,\,v\to\infty\,\,\,\,{\rm for}\,\,k \geq 1.  
$$
Moreover, we have 
\begin{equation}
\frac{\phi^{[\alpha_1]}(v)\cdots \phi^{[\alpha_{\mu}]}(v)}
{\phi(v)^{\mu}} 
\sim 
v^{\frac{(2m+1)\mu}{2m-1}} \cdot
\sum_{j=0}^{\infty}c_j v^{-\frac{2m}{2m-1}j}
\,\,\,\,\,\,\,{\rm as}\,\,v\to\infty.  
\label{eqn:6.162}
\end{equation}
Therefore (\ref{eqn:6.16}),(\ref{eqn:6.162}) and Lemma 6.1 imply 
$$
a_{\mu}(v) 
\sim 
v^{\frac{m-1+(2m+1)\mu}{2m-1}} \cdot 
\exp\{-a v^{\frac{2m}{2m-1}}\}
\cdot\sum_{j=0}^{\infty}c_j v^{-\frac{2m}{2m-1}j}
\,\,\,\,\,\,\,{\rm as}\,\,v\to\infty. 
$$

%\vspace{.5em}
%%%%%%%%%%%%

{\it Proof of Lemma 6.1}.\,\,\,
Changing the integral variable, we have 
\begin{equation}
\phi_l(v)=v^{\frac{l+1}{2m-1}} \int_{-\infty}^{\infty} 
t^l e^{-\tilde{v}p(t)}dt, 
\label{eqn:6.163}
\end{equation}
where $\tilde{v}=v^{\frac{2m}{2m-1}}$ and 
$p(t)=t^{2m}-t$. 
We divide (\ref{eqn:6.163}) into two parts: 
\begin{equation}
\phi_l(v)=v^{\frac{l+1}{2m-1}}
\{I_1(\tilde{v})+I_2(\tilde{v})\},
\label{eqn:6.17}
\end{equation}
with 
$$
I_1(\tilde{v})
= 
\int_{|t-\alpha|\leq\delta} 
t^l e^{-\tilde{v}p(t)}dt \,\,\, \mbox{and} \,\,\,
I_2(\tilde{v})
= 
\int_{|t-\alpha|>\delta} 
t^l e^{-\tilde{v}p(t)}dt,
$$
where  $\delta>0$ is small and 
$\alpha=(2m)^{\frac{-1}{2m-1}}$. 
Note that $p'(\alpha)=0$. 

%\vspace{.5em}
First we consider the function $I_1$. 
By Taylor's formula, we have   
\begin{eqnarray*}
I_1(\tilde{v})
\!\!\!&=&\!\!\!
\int_{|t-\alpha|<\delta} 
t^l \exp\{
-\tilde{v}[p(\alpha)+\tilde{p}(t-\alpha)(t-\alpha)^2]
\}dt  \\ 
 \!\!\!&=&\!\!\!
e^{a\tilde{v}}\int_{|t|\leq\delta} 
(t+\alpha)^l e^{-\tilde{v}\tilde{p}(t)t^2}dt,  
\end{eqnarray*}
where 
$a=-p(\alpha)
=[(2m)^{\frac{-1}{2m-1}}-(2m)^{\frac{-2m}{2m-1}}]>0$ and 
$\tilde{p}(t)=\int_0^1 (1-u)p''(ut+\alpha) du$. 
Set $s=\tilde{p}(t)^{\frac{1}{2}}t$ for $|t|\leq\delta$ and 
$\delta_{\pm}=\tilde{p}(\pm \delta)^{\frac{1}{2}}\delta$, 
respectively. 
Then there is a function 
$\hat{p}\in C^{\infty}([-\delta_-,\delta_+])$ such that 
$t=\hat{p}(s)$ and $\hat{p}'>0$. 
Changing the integral variable, 
we have   
$$
I_1(\tilde{v})=\int_{-\delta_-}^{\delta_+} 
e^{-\tilde{v}s^2} \Psi(s) ds, 
$$
where $\Psi(s)=(\hat{p}(s)+\alpha)^l \cdot \hat{p}'(s)$.  
Since $\Psi\in C^{\infty}([-\delta_-,\delta_+])$, we have 
\begin{eqnarray}
I_1(v)
\!\!\!&=&\!\!\!
\tilde{v}^{-\frac{1}{2}}
\int_{-\delta_- \tilde{v}^{\frac{1}{2}}}
^{\delta_+ \tilde{v}^{\frac{1}{2}}} 
e^{-u^2} \Psi(u \tilde{v}^{-\frac{1}{2}})du \nonumber\\
\!\!\!&\sim&\!\!\!
\tilde{v}^{-\frac{1}{2}}
\sum_{j=0}^{\infty} c_j \tilde{v}^{-j} \quad 
\mbox{as $\tilde{v}\to\infty$.}
\label{eqn:6.22}
\end{eqnarray}
We remark that $\int_{-\infty}^{\infty}e^{-u^2}u^j du=0$ if 
$j \in {\bf Z}$ is odd. 

%\vspace{.5em} 
Next we consider the function $I_2$. 
Let $p_d$ be the function defined by 
$      
p_d(t)= d|t-\alpha|-a
$
where $d>0$. 
We can choose $d>0$ such that 
$p(t) \geq p_d(t)$ for $|t-\alpha|>\delta$. 
Then we have 
\begin{eqnarray}
|I_2(\tilde{v})| 
&\leq& 
\int_{|t-\alpha|>\delta} e^{-\tilde{v}p_d(t)} dt \nonumber\\
&\leq& 
2 C e^{a\tilde{v}}
\int_{\delta}^{\infty} e^{-d\tilde{v}t} dt 
\leq 
2C\tilde{v}^{-1} e^{[a-d\delta]\tilde{v}}.  
\label{eqn:6.23}
\end{eqnarray}

Finally putting 
(\ref{eqn:6.17}),(\ref{eqn:6.22}),(\ref{eqn:6.23}) together, 
we have the asymptotic expansion in Lemma 6.1. 
\qed

%%%%%%%%%%%%%%%%%%%%%%
%%%%%%%%%%%%%%%%%%%%%%
%%%%%%%%%%%%%%%%%%%%%%
%%%%%%%%%%%%%%%%%%%%%%
%%%%%%%%%%%%%%%%%%%%%%

\subsection{Asymptotic expansion of $L_{\mu}$}

Let $A \in C^{\infty}({\bf R})$ be an even or odd  function 
(i.e.  $A(-v)=A(v)$ or $-A(v)$) satisfying  
$$
A(v)\sim v^{\frac{n}{2m-1}}\cdot\exp\{-av^{\frac{2m}{2m-1}}\}
\cdot\sum_{j=0}^{\infty}c_j v^{-\frac{2m}{2m-1}j} 
\quad \mbox{as $v\to +\infty$,}
$$
where $n\in {\bf N}$ and the constant $a$ is as in Lemma 6.1. 
Let $L\in C^{\omega}({\bf R})$ be the function defined by 
\begin{equation}
L(u)=\int_{-\infty}^{\infty} e^{uv} A(v) dv.
\label{eqn:6.230}
\end{equation}
In this section we give the asymptotic expansion of $L$ at infinity. 
%%%%%%%%%%%%%%%%%%%%%%%%%%%%%%%%%%%%%%
\begin{lemma} \quad\,
$
L(u)
\sim u^{m-1+n} \cdot e^{u^{2m}}\cdot 
\sum_{j=0}^{\infty} c_j u^{-2mj} \,\,\,\,\,\,\,{\rm as}\,\,
u \to +\infty. 
$
\end{lemma}
%%%%%%%%%%%%%%%%%%%%%%%%%%%%%%%%%%%%%%%

{\it Remark}.\,\,\, 
Lemmas 6.1, 6.2 imply that for $\mu,l \geq 0$, 
\begin{eqnarray}
&&L_{\mu}^{(l)}(u)
 =  
\frac{d^l}{du^l}L_{\mu}(u) \nonumber\\
&&\quad\sim
u^{2m-2+(2m+1)\mu+(2m-1)l} \cdot e^{u^{2m}}\cdot 
\sum_{j=0}^{\infty} c_j u^{-2mj} \quad \mbox{ as 
$u \to \infty$.}
\label{eqn:6.24}  
\end{eqnarray} 

{\it Proof}. \,\,\,
We only show Lemma 6.2 in the case where 
$A$ is an even function. 
Let $\hat{A}\in C^{\infty}({\bf R}\setminus\{0\})$ be defined by 
\begin{equation}
A(v)=v^{\frac{n}{2m-1}}\cdot\exp\{-a v^{\frac{2m}{2m-1}}\}
\cdot \hat{A}(v^{\frac{2m}{2m-1}}), 
\label{eqn:6.241}
\end{equation}
then $\hat{A}(x)\sim\sum_{j=0}^{\infty}c_j x^{-j}$ as $x\to\infty$.  
Substituting (\ref{eqn:6.241}) into (\ref{eqn:6.230}), we have 
$$
L(u)=\int_{-\infty}^{\infty} \exp\{-a|v|^{\frac{2m}{2m-1}} 
+uv\} |v|^{\frac{n}{2m-1}} \hat{A}(|v|^{\frac{2m}{2m-1}}) dv. 
$$
Changing the integral variable and setting 
$q(t)=at^{2m}-t^{2m-1}$, 
we have 
\begin{equation}
L(u)=
(2m-1)u^{2m+n-1} \int_{-\infty}^{\infty}
e^{-u^{2m}q(t)} \hat{A}(t^{2m}u^{2m}) t^{2m+n-2}dt.
\label{eqn:6.240}
\end{equation}
Now we divide the integral in (\ref{eqn:6.240}) 
into two parts: 
\begin{equation}
L(u)=(2m-1)u^{2m+n-1}\{
J_1(\tilde{u})+J_2(\tilde{u})  
\}, 
\label{eqn:6.25} 
\end{equation}
with 
\begin{eqnarray}
J_1(\tilde{u})\!\!\!&=&\!\!\! \int_{|t-\beta|\leq\delta} e^{-\tilde{u}q(t)} 
               \hat{A}(\tilde{u}t^{2m})t^{2m+n-2}dt,
               \nonumber\\  
J_2(\tilde{u})\!\!\!&=&\!\!\! \int_{|t-\beta|>\delta} e^{-\tilde{u}q(t)} 
               \hat{A}(\tilde{u}t^{2m})t^{2m+n-2}dt,  
               \nonumber 
\end{eqnarray}
where $\tilde{u}=u^{2m}$, $\delta>0$ is small and 
$\beta=(2m\!-\!1)\cdot(2m a)^{-1}$. Note that $q'(\beta)=0$. 

%\vspace{.3em}
First we consider the function $J_1$. 
By Taylor's formula, we have   
$$
J_1(\tilde{u})=e^{\tilde{u}}\int_{|t|\leq\delta}
e^{-\tilde{u}\tilde{q}(t)t^2} 
\hat{A}(\tilde{u}(t+\beta)^{2m})(t+\beta)^{2m+n-2}dt, 
$$
where $\tilde{q}(t)=\int_0^1(1-v)q''(vt+\alpha)dv$. 
Note that $q(0)=-1$. 
Set $s=\tilde{q}(t)^{\frac{1}{2}}t$ for $|t|\leq\delta$ 
and $\tilde{\delta}_{\pm}=\tilde{q}(\pm\delta)^{\frac{1}{2}}\delta$, 
respectively. 
Then there is a function 
$\hat{q}\in C^{\infty}([-\tilde{\delta}_-,\tilde{\delta}_+])$ 
such that $t=\hat{q}(s)$ and $\hat{q}'>0$. 
Changing the integral variable, we have 
\begin{equation}
J_1(\tilde{u})=\int_{-\tilde{\delta}_-}^{\tilde{\delta}_+} 
e^{-\tilde{u}s^2} \tilde{\Psi}(s,\tilde{u}) ds, 
\label{eqn:6.280}
\end{equation}
where  
$
\tilde{\Psi}(s,\tilde{u})=\hat{A}(\tilde{u}(\hat{q}(s)+\beta)^{2m}) 
(\hat{q}(s)+\beta)^{2m+n-2} \hat{q}'(s). 
$
Since $\hat{A}(x)\sim \sum_{j=0}^{\infty}c_j x^{-j}$ as $x\to \infty$, 
we have 
\begin{equation}
\tilde{\Psi}(s,\tilde{u})\sim 
\sum_{j=0}^{\infty}c_j(s)\tilde{u}^{-j}\quad \mbox{as  
$\tilde{u}\to\infty$}, 
\label{eqn:6.29}
\end{equation} 
where 
$c_j \in C^{\infty}([-\tilde{\delta}_-,\tilde{\delta}_+])$. 
Substituting (\ref{eqn:6.29}) into (\ref{eqn:6.280}), 
we have 
\begin{equation}
J_1(\tilde{u})\sim \tilde{u}^{-\frac{1}{2}}
e^{\tilde{u}}\sum_{j=0}^{\infty} c_j \tilde{u}^{-j} 
\quad \mbox{as $\tilde{u}\to\infty$}. 
\label{eqn:6.30} 
\end{equation}

Next we consider the function $J_2$. 
By a similar argument about 
the estimate of $I_2(\tilde{v})$ in the proof of Lemma 6.1, 
we can obtain 
\begin{equation}
|J_2(\tilde{u})|\leq C \tilde{u}^{-1} 
e^{[1-\varepsilon]\tilde{u}}, 
\label{eqn:6.31} 
\end{equation}
where $\varepsilon$ is a positive constant. 

Finally putting 
(\ref{eqn:6.25}),(\ref{eqn:6.30}),(\ref{eqn:6.31}) together, 
we obtain the asymptotic expansion in Lemma 6.2. 
\qed

%%%%%%%%%%%%%%%%%%%%%%%
%%%%%%%%%%%%%%%%%%%%%%%

\subsection{Proof of Proposition 6.1} 

We can construct the function $h\in C^{\infty}([0,\infty))$ 
such that if $Y=\tilde{f}(X)^{\frac{1}{2m}}$, 
then $X=Yh(Y)$. 
In fact $\frac{d}{dX}[\tilde{f}(X)^{\frac{1}{2m}}]>0$ for 
$X \geq 0$. 
Set $t=\tilde{f}(x)^{\frac{1}{2m}}\xi^{-1}$.   
Then we can write $t_0=th(t \xi)$. 
Note that $\hat{g}(X)^{\frac{1}{2m}} \cdot h(X)=1$ 
for $X, Y \geq 0$ 
and hence $h'(Y) \geq 0$ for $Y \geq 0$. 
Let us prepare two lemmas 
for the proof of Propositions 6.1, 6.2. 

%%%%%%%%%%%%%%%%%%%%%%%%%%%%%%%%%%%%%%%%%%%%%%%%%%%%%%
\begin{lemma}
\,\,  Assume that the functions $a,b$ and $c$ 
on $[0,\varepsilon) \times [0,1]$ 
satisfy $a(\xi,t_0)=b(\xi,t)=c(\xi,\tau)$.   
If one of these functions 
belongs to $C^{\infty}([0,\varepsilon)\times [0,1])$, 
then so do the others.  
\end{lemma} 
%%%%%%%%%%%%%%%%%%%%%%%%%%%%%%%%%%%%%%%%%%%%%%%%%%%%%%%

{\it Proof}.\,\,\, 
This lemma is directly shown by 
the relation between three variables $t_0,t$ and $\tau$.  
\qed

%%%%%%%%%%%%%%%%%%%%%%%%%%%%%%%%%%%%%%%%%%%%%%%%%%%%%
\begin{lemma}
There is a positive number $\alpha$ such that 
$1-t_0^{2m} \geq \tau-\alpha\xi$. 
\end{lemma}
%%%%%%%%%%%%%%%%%%%%%%%%%%%%%%%%%%%%%%%%%%%%%%%%%%%%%
We remark that the above constant  $\alpha$ is same as  
that in Proposition 6.2. 

{\it Proof}.\,\,\, 
By definition, we have 
$$
1-t_0^{2m}=1-t^{2m}h(t \xi)^{2m}=(1-t^{2m})-(h(t\xi)-1)t^{2m}.$$
Since $h(0)=1$ and $h'(X)\geq 0$, 
we have 
$
h(t\xi)^{2m}-1\leq \alpha t\xi
$
for some positive number $\alpha$. 
Therefore we have 
$$
1-t_0^{2m}\geq \tau-\alpha t^{2m+1}\xi 
\geq \tau-\alpha\xi.
$$
\qed
%%%%%

{\it Proof of Proposition 6.1.} 
\,(i) \, Recall the definition of the function 
$K_{\mu}^{\langle 1 \rangle}$: 
\begin{equation}
K_{\mu}^{\langle 1 \rangle}(\tau,\xi)
=
\int_{1}^{\infty}
e^{-s^{2m}}L_{\mu}(t_0 s)s^{4m+1-\mu}ds,
\label{eqn:6.32}
\end{equation}
where 
$$
L_{\mu}(u) 
=
\int_{-\infty}^{\infty}
e^{uv}a_{\mu}(v) dv \quad {\rm and} \quad
a_{\mu}(v) =  
\left. \frac{1}{\mu!}
\frac{\partial^{\mu}}{\partial X^{\mu}} 
\frac{1}{\phi(v,X)}\right|_{X=0}. 
$$
%%%%%%%%%%%%%%%%%%%%%%%%%%%%%%%%%%%%%%%

We obtain 
the Taylor expansion of $L_{\mu}(t_0;s)=L_{\mu}(tsh(t\xi))$ 
with respect to $\xi$: 
\begin{equation}
L_{\mu}(t_0 s)=
\sum_{k=0}^{k_0} L_{\mu,k}(t;s) \xi^k 
+\tilde{L}_{\mu,k_0}(t,\xi;s)\xi^{k_0+1}, 
\label{eqn:6.330}
\end{equation}
where 
\begin{eqnarray}
L_{\mu,k}(t;s)
\!\!\!&=&\!\!\!
\frac{1}{k!} 
\left.\frac{\partial^k}{\partial \xi^k}
 L_{\mu}(ts h(t \xi))\right|_{\xi=0}, 
\label{eqn:6.34}\\
\tilde{L}_{\mu,k_0}(t,\xi;s)
\!\!\!&=&\!\!\!
\frac{1}{k_0!} 
\int_0^1 (1-p)^{k_0} \left.
\frac{\partial^{k_0+1}}{\partial X^{k_0+1}}
 L_{\mu}(ts h(t X))\right|_{X=\xi p}   dp. 
\label{eqn:6.340} 
\end{eqnarray}

Substituting (\ref{eqn:6.330}) into (\ref{eqn:6.32}), 
we have 
$$
K_{\mu}^{\langle 1 \rangle}(\tau,\xi) 
=  
\sum_{k=0}^{k_0} K_{\mu,k}(t) \xi^{k}
+\tilde{K}_{\mu,k_0}(t,\xi)\xi^{k_0+1}, 
$$
where 
\begin{eqnarray}
K_{\mu,k}(t) 
\!\!\!&=&\!\!\!
\int_1^{\infty} e^{-s^{2m}} L_{\mu,k}(t;s)s^{4m+1-\mu}ds
\label{eqn:6.37}
\\
\tilde{K}_{\mu,k_0}(t,\xi) 
\!\!\!&=&\!\!\! 
\int_{1}^{\infty} e^{-s^{2m}}\tilde{L}_{\mu,k_0}
(t,\xi;s)s^{4m+1-\mu}ds. 
\label{eqn:6.38} 
\end{eqnarray}

%\vspace{.5em}

First we consider the singularity of $K_{\mu,k}$ at $t=1$.

By a direct computation in (\ref{eqn:6.34}), 
we have   
$$
L_{\mu,k}(t;s)=
\sum_{l=1}^{k} h_{l}(t) s^l L_{\mu}^{(l)}(ts), 
$$
where $h_{l} \in C^{\infty}([0,1])$ 
(which depends on $\mu,k$). 
We define the function $S_{\mu,k}(t;s^{2m})$ by 
\begin{equation}
L_{\mu,k}(t;s)=
s^{2m-2+(2m+1)\mu +2mk}\cdot e^{t^{2m}s^{2m}}\cdot 
S_{\mu,k}(t;s^{2m}).
\label{eqn:6.380}
\end{equation}
Lemma 6.2 implies  
\begin{equation} 
S_{\mu,k}(t;s^{2m}) \sim \sum_{j=0}^{\infty} c_j(t) s^{-2mj} 
\,\,\,\,\,\,{\rm as}\,\,\,s \to \infty, 
\label{eqn:6.381} 
\end{equation}
where $c_j \in C^{\infty}([0,1])$. 
 
Substituting (\ref{eqn:6.380}) into (\ref{eqn:6.37}), 
we have 
\begin{eqnarray}
K_{\mu,k}(t) 
\!\!\!&=&\!\!\! 
\int_{1}^{\infty} e^{-[1-t^{2m}]s^{2m}} 
S_{\mu,k}(t;s^{2m})s^{2m-2+2m\mu+2mk} ds 
\nonumber\\
\!\!\!&=&\!\!\!
\frac{1}{2m} \int_{1}^{\infty} 
e^{-\chi^{-1}(\tau) \sigma} 
S_{\mu,k}(t;\sigma)\sigma^{2+\mu+k} d\sigma. 
\label{eqn:6.40}
\end{eqnarray}
Moreover substituting (\ref{eqn:6.381}) into (\ref{eqn:6.40}), 
we have 
$$
K_{\mu,k}(t)=\frac{\varphi_{\mu,k}(\tau)}{\tau^{3+\mu+k}} 
+\psi_{\mu,k}(\tau) \log \tau, 
$$
where $\varphi_{\mu,k},\psi_{\mu,k} \in C^{\infty}([0,1])$.

%\vspace{.5em}
Next we obtain the inequality 
$|\tilde{K}_{\mu,k}^{\langle 1 \rangle}(\tau,\xi)| \leq 
C_{\mu,k_0}[\tau-\alpha\xi]^{-4-\mu-k_0}$ for 
some positive constant $C_{\mu,k_0}$. 
By a direct computation in (\ref{eqn:6.340}), 
we have
$$
\frac{\partial^{k_0+1}}{\partial X^{k_0+1}}
L_{\mu}(tsh(tX))=
\sum_{l=1}^{k_0+1} 
\tilde{h}_l (t,X) s^l L_{\mu}^{(l)}(tsh(tX)), 
$$
where $\tilde{h}_l$ are bounded functions 
(depending on $\mu,k_0$). 
Since $h'(X) \geq0$,  
we can obtain 
\begin{eqnarray}
|\tilde{L}_{\mu,k_0}(t,\xi;s)| 
\!\!\!&\leq&\!\!\! 
C \sum_{l=1}^{k_0+1} s^{l}L_{\mu}^{(l)}(tsh(ts))\nonumber\\
\!\!\!&\leq&\!\!\! 
C s^{k_0+1} L_{\mu}^{(k_0+1)}(t_0 s)\nonumber\\
\!\!\!&\leq&\!\!\!
C s^{2m-2+(2m+1)\mu+2m(k_0+1)}\cdot e^{t_0^{2m}s^{2m}}
\label{eqn:6.50}
\end{eqnarray}
for $s \geq 1$.
Substituting (\ref{eqn:6.50}) to (\ref{eqn:6.38}), 
we obtain 
\begin{eqnarray*}
\tilde{K}_{\mu,k_0}(t,\xi) 
\!\!\!&\leq&\!\!\! 
C\int_{1}^{\infty}
e^{-[1-t_0^{2m}]s^{2m}} s^{6m-2+2m\mu+2m(k_0+1)} ds \\
\!\!\!&\leq&\!\!\!
C [1-t_0^{2m}]^{-4-\mu-k_0} 
 \leq 
C [\tau-\alpha\xi]^{-4-\mu-k_0} 
\end{eqnarray*}
by Lemma 6.4. 
This completes the proof of Proposition 6.2 (i). 

%%%%%%%%%%%%%%%%%%%%%%%%%%%%%%%%%%%%%%%%%%%%
%\vspace{.5em}

(ii)\,\,
Recall the definition of the function  
$\tilde{K}_{\mu}^{\langle 1 \rangle}$:  
\begin{equation}
\tilde{K}_{\mu_0}^{\langle 1 \rangle}(\tau,\xi)
=
\int_{1}^{\infty}
e^{-s^{2m}} \tilde{L}_{\mu_0}(t_0,\xi;s)s^{4m-\mu_0}ds,
\label{eqn:6.52}
\end{equation}
where
\begin{eqnarray}
&&\tilde{L}_{\mu_0}(t_0,\xi;s)
=
\int_{-\infty}^{\infty} e^{t_0 s v} r_{\mu_0}(v,\xi s^{-1})dv, 
\label{eqn:6.53}\\
&&r_{\mu_0}(v,X) 
= 
\frac{1}{\mu_0 !} 
\int_0^1(1-p)^{\mu_0} \left.
\frac{\partial^{\mu_0}}{\partial Y^{\mu_0+1}} 
\frac{1}{\phi(v,Y)}\right|_{Y=Xp} dp. 
\label{eqn:6.54}
\end{eqnarray}

The following lemma is necessary to 
obtain the estimate of $\tilde{K}_{\mu}^{\langle 1 \rangle}$ 
in (ii). 
%%%%%%%%%%%%%%%%%%%%%%%%%%%
\begin{lemma}
\,\,\,
$
|r_{\mu_0}(v,X)|< C |v|^{\frac{(2m+1)(\mu_0+1)+m-1}{2m-1}} 
e^{-a |v|^{\frac{2m}{2m-1}}}$\,\, for $v\in {\bf R}$,  $X\geq0$. 
\end{lemma}
%%%%%%%%%%%%%%%%%%%%%%%%%%%
We remark that the constant $a$ is as in Lemma 6.1. 
The proof of the above lemma is given soon later. 

Applying Lemma 6.5 to (\ref{eqn:6.53}),  
we have 
\begin{eqnarray}
|\tilde{L}_{\mu_0}(t_0,\xi,s)| 
\!\!\!&\leq&\!\!\! 
C \int_{-\infty}^{\infty} 
|v|^{\frac{(2m+1)(\mu_0+1)+m-1}{2m-1}} 
e^{t_0 sv-a |v|^{\frac{2m}{2m-1}}}dv \nonumber\\
\!\!\!&\leq&\!\!\! 
C s^{2m-2+(2m+1)(\mu_0+1)} e^{t_0^{2m} s^{2m}}. 
\label{eqn:6.55}
\end{eqnarray}
The second inequality is given by Lemma 6.2. 
Moreover substituting (\ref{eqn:6.55}) into (\ref{eqn:6.52}), 
we have 
\begin{eqnarray*}
|\tilde{K}_{\mu_0}^{\langle 1 \rangle}(\tau,\xi)| 
\!\!\!&\leq&\!\!\!  C \int_1^{\infty} 
         e^{-[1-t_0^{2m}]s^{2m}}s^{8m-1+2m\mu_0} ds \\
\!\!\!&\leq\!\!\!&  C [1-t_0^{2m}]^{-4-\mu_0} 
\leq C [\tau-\alpha \xi]^{-4-\mu_0}, 
\end{eqnarray*}
by Lemma 6.4. 
Therefore  
we obtain the estimate of $\tilde{K}_{\mu_0}^{\langle 1 \rangle}$ 
in Proposition 6.1 (ii). 
\qed

{\it Proof of Lemma 6.5}.\,\, 
We only consider the case where $v$ is positive. 
The proof for the case where $v$ is negative 
is given in the same way.

By a direct computation, we have 
\begin{equation}
\frac{\partial^{\mu_0+1}}{\partial X^{\mu_0+1}}
\frac{1}{\phi(v,X)}
=  
\sum_{|\alpha|=\mu_0+1}
C_{\alpha} \frac{
\phi^{[\alpha_1]}(v,X)\cdots \phi^{[\alpha_{\mu_0+1}]}(v,X)
}{\phi(v,X)^{\mu_0+2}}, 
\label{eqn:6.550} 
\end{equation} 
where $\phi^{[k]}(v,X)= 
\frac{\partial^k}{\partial X^k}\phi(v,X)$ and  
$C_{\alpha}$'s are constants depending on 
$\alpha=(\alpha_1,\ldots,\alpha_{\mu_0+1})
\in {\bf Z}_{\geq 0}^{\mu_0+1}$. 
By a direct computation, 
the function $\phi^{[k]}$'s are  expressed in the form: 
\begin{equation}
\phi^{[k]}(v,X)=
\sum_{|\beta|=k} c_{\beta} F_{\beta}(\tilde{v},\tilde{X}),
\label{eqn:6.56}
\end{equation}
where $\tilde{v}=v^{\frac{2m}{2m-1}}$, 
$\tilde{X}=X v^{\frac{-1}{2m}}$, 
$c_{\beta}$'s are constants and 
\begin{eqnarray}
F_{\beta}(\tilde{v},\tilde{X})
\!\!\!&=&\!\!\!
\int_{-\infty}^{\infty}
e^{-\tilde{g}(X w)w^{2m}+vw}
w^{\gamma} \prod_{\beta} 
\hat{g}^{(\beta_k)}(Xw) dw \nonumber\\
\!\!\!&=&\!\!\!
\tilde{v}^{\frac{\gamma+1}{2m}}
\int_{-\infty}^{\infty}
e^{-\tilde{v} p(s,\tilde{X})}
s^{\gamma}\prod_{\beta} 
{\hat{g}}^{(\beta_k)}(\tilde{X}s)  ds, 
\label{eqn:6.57}
\end{eqnarray}
with $p(s,\tilde{X}) =\hat{g}(\tilde{X}s)s^{2m}-s$, 
$\beta_k \in {\bf N}$ 
and $\gamma \in {\bf N}$ depending on $\beta=(\beta_k)_k$.  
In order to apply the stationary phase method to the above integral, 
we must know the location of the critical points of 
the function $p(\cdot,\tilde{X})$. 
The lemma below gives the information about it. 
%%%%%%%%%%%%%%%%%%%%%%%%%
\begin{lemma}
There exists a function $\alpha \in C^{\infty}([0,\infty))$ 
such that 
\begin{equation}
\frac{\partial}{\partial s} p(\alpha(\tilde{X}), \tilde{X})=0 \, 
\mbox{ and 
$\alpha_0 \leq \alpha(\tilde{X}) \le \alpha_1$  \,
for \, $\tilde{X} \geq 0$}, 
\label{eqn:6.58}
\end{equation}
where $\alpha_0=(2m)^{\frac{-1}{2m-1}}$ and 
 $\alpha_1=(\frac{4}{5}2m)^{\frac{-1}{2m-1}}$. 
\end{lemma}
%%%%%%%%%%%%%%%%%%%%%%%%%%%%

{\it Proof}. \,\,\, 
By a direct computation, 
we have 
\begin{eqnarray*}
\frac{\partial}{\partial s} p(s,\tilde{X}) 
\!\!\!&=&\!\!\!
s^{2m-1}[\hat{g}'(\eta)\eta+ 2m \hat{g}(\eta)]-1,\\
\frac{\partial^2}{\partial s^2} p(s,\tilde{X}) 
\!\!\!&=&\!\!\!
s^{2m-2}[\hat{g}''(\eta)\eta^2+ 4m \hat{g}'(\eta)\eta 
+2m(2m-1)\hat{g}(\eta)],
\end{eqnarray*}
where $\eta=\tilde{X}s$. 
It is easy to
obtain the following inequalities 
by using the conditions (\ref{eqn:6.1}),(\ref{eqn:6.2}). 
\begin{eqnarray}
&&\frac{\partial}{\partial s} p(\alpha_0,\tilde{X}) \leq 0 < 
 \frac{\partial}{\partial s} p(\alpha_1,\tilde{X}) 
\quad \mbox{for $\tilde{X}\geq 0$},
\label{eqn:6.61}\\
&&\frac{\partial^2}{\partial s^2} p(s,\tilde{X})\geq c >0 
\quad \mbox{on $[\alpha_0,\alpha_1]\times [0,\infty)$}.
\label{eqn:6.62}
\end{eqnarray}
for some positive constant $c$. 
Then the inequalities (\ref{eqn:6.61}),(\ref{eqn:6.62}) 
imply the claim in Lemma 6.6  
by the implicit function theorem. 
\qed
%%% 

Now we divide the integral in (\ref{eqn:6.57}) into two parts. 
$$
F_{\beta}(\tilde{v},\tilde{X}) 
= \tilde{v}^{\frac{\gamma+1}{2m}} \{
I_1(\tilde{v}, \tilde{X})+I_2(\tilde{v}, \tilde{X}) 
\}, 
$$
where 
\begin{eqnarray*}
I_1(\tilde{v}, \tilde{X})
\!\!\!&=&\!\!\!
\int_{|t-\alpha(\tilde{X})| \leq \delta} 
e^{-\tilde{v}p(t,\tilde{X})t^2} 
t^{\gamma}
\prod_{\beta} 
{\hat{g}}^{(\beta_k)}(\tilde{X}t) 
dt, \\
I_2(\tilde{v}, \tilde{X})
\!\!\!&=&\!\!\!
\int_{|t-\alpha(\tilde{X})| > \delta} 
e^{-\tilde{v}p(t,\tilde{X})t^2} 
t^{\gamma}
\prod_{\beta} 
{\hat{g}}^{(\beta_k)}(\tilde{X}t) 
dt, 
\end{eqnarray*}
where $\delta>0$ is small. 

%\vspace{.5em} 
First we consider the function $I_1$. 
By Lemma 6.6 and Taylor's formula, we have 
\begin{equation}
I_1(\tilde{v}, \tilde{X})
=    
e^{a(\tilde{X})\tilde{v}} 
\int_{|t|\leq \delta}
e^{-\tilde{v}\tilde{p}(t,\tilde{X})t^2} 
(t+\alpha(\tilde{X}))^{\gamma}
\prod_{\beta}\hat{g}^{(\beta_k)}
(\tilde{X}(t+\alpha(\tilde{X}))) 
dt,  
\end{equation}
where $a(\tilde{X})=-p(\alpha(\tilde{X}),\tilde{X})$ and 
$\tilde{p}(t,\tilde{X}) =
   \int_{0}^{1}(1-u)
   \frac{\partial^2}{\partial s^2}p(ut+\alpha(\tilde{X}),\tilde{X}) du$. 
Set $s=\tilde{p}(t,\tilde{X})^{\frac{1}{2}} t$ on 
$[-\delta,\delta]\times [0,\infty)$ and 
$\delta_{\pm}(\tilde{X})= 
\tilde{p}(\pm\delta, \tilde{X})^{\frac{1}{2}}\delta$, 
respectively. 
Then there is a function $\hat{p}\in C^{\infty}
([-\delta_-(\tilde{X}),\delta_+(\tilde{X})]\times [0,\infty))$ 
such that 
$t=\hat{p}(s,\tilde{X})$ and 
$\frac{\partial}{\partial s} \hat{p}(s,\tilde{X})>0$.  
Changing the integral variable, 
we have
$$
I_1(\tilde{v}, \tilde{X})
=   
e^{a(\tilde{X})\tilde{v}} 
\int_{-\delta_-(\tilde{X})}^{\delta_+(\tilde{X})}
e^{-\tilde{v}s^2}
\Psi(s,\tilde{X})ds, 
$$
where $\Psi(s,\tilde{X})
=       
(\hat{p}(s,\tilde{X})+\alpha(\tilde{X}))^{\gamma}\cdot
\prod_{\beta}\hat{g}^{(\beta_k)}
(\tilde{X}(\hat{p}(s,\tilde{X})+\alpha(\tilde{X})))\cdot  
\frac{\partial}{\partial s}\hat{p}(s,\tilde{X}). 
$      
Since $\Psi \in C^{\infty}
([-\delta_-(\tilde{X}),\delta_+(\tilde{X})]\times [0,\infty))$, 
we have 
\begin{eqnarray}
I_1(\tilde{v}, \tilde{X}) \cdot 
\{\tilde{v}^{-\frac{1}{2}}
e^{a(\tilde{X})\tilde{v}}\}^{-1}
\!\!\!&=&\!\!\!
\int_{-\delta_-(\tilde{X})\tilde{v}^{\frac{1}{2}}}
^{\delta_+(\tilde{X})\tilde{v}^{\frac{1}{2}}}
e^{-u^2}
\Psi(u \tilde{v}^{-\frac{1}{2}}, \tilde{X})du \nonumber\\    
\!\!\!&\to&\!\!\! 
\sqrt{\pi}\Psi(0,\tilde{X}) \quad {\rm as}\,\, \tilde{v}\to\infty. 
\label{eqn:6.68}
\end{eqnarray}
Note that 
$\Psi(0,\tilde{X}) =\alpha(\tilde{X})^{\gamma}\cdot
\{\frac{1}{2} 
\frac{\partial^2}{\partial t^2} 
p(\alpha(\tilde{X}),\tilde{X})\}^{-\frac{1}{2}} \cdot
\prod_{\beta}\hat{g}
^{(\beta_k)}(\tilde{X} \alpha(\tilde{X}))$.

%\vspace{.5em}
Next we consider the function $I_2$. 
In a similar argument about the estimate of $I_2(\tilde{v})$ 
in the proof in Lemma 6.1, 
we can obtain 
\begin{equation}
|I_2(\tilde{v},\tilde{X})| 
\leq C \tilde{v}^{-1} e^{a(\tilde{X})\tilde{v}
-\varepsilon\tilde{v}}
\label{eqn:6.69}
\end{equation}
where $\varepsilon$ is a positive constant. 

%%%%%%%%%%%%%%%%%%%%%%%%%%%%%%%%%
%\begin{lemma}
%There is a positive constant $b$ such that 
%$\tilde{P}(s,\tilde{X})\geq b$. 
%\end{lemma}
%%%%%%%%%%%%%%%%%%%%%%%%%%%%%%%%
%
%{\it Proof}.\,\,\, 
%Since $\frac{\partial^2}{\partial s^2} 
%p(us+ \alpha(\tilde{X}),\tilde{X})\geq 
%c(us+\alpha(\tilde{X}))^{2m-2}$ for some $c>0$ by (), 
%we have 
%\begin{eqnarray}
%\tilde{P}(s,\tilde{X})
%\!\!\!&\geq&\!\!\! 
%c \int_0^1 (1-u)(su-\alpha(\tilde{X}))^{2m-2} du\\
%\!\!\!&\geq&\!\!\!
%c \left\{
%\int_0^1 (1-u)^{\frac{1}{m-1}}(su-\alpha(\tilde{X}))^{2} du 
%\right\}^{m-1}\\
%\!\!\!&\geq&\!\!\!
%c \left\{
%\frac{(s+2\alpha(\tilde{X}))^2}{12}+ 
%\frac{\alpha(\tilde{X})^2}{6} 
%\right\}^{m-1} \geq c \frac{\alpha(\tilde{X})^{2m-2}}{6^{m-1}}>0 
%\end{eqnarray}
%by H\"older's inequality. 
%If we set $b=c\, 6^{1-m} \alpha_0^{2m-2}$, 
%we  obtain the lemma. 
%\qed

Putting (\ref{eqn:6.68}),(\ref{eqn:6.69}) together, 
we have  
\begin{equation}
\lim_{\tilde{v} \to \infty} 
F_{\beta}(\tilde{v},\tilde{X})\cdot 
\{ 
\tilde{v}^{\frac{\gamma+1-m}{2m}} 
e^{-p(\tilde{X})\tilde{v}}
\}^{-1}=
\sqrt{\pi}  \tilde{\Psi}_{\beta}(0,\tilde{X}). 
\label{eqn:6.70}
\end{equation}

Now under the condition $|\beta|=k$, 
the number $\gamma$ in (\ref{eqn:6.57}) attains the maximum value 
$(2m+1)k$ when $\beta=(1,\ldots,1)$. 
Therefore (\ref{eqn:6.56}),(\ref{eqn:6.70}) imply that 
\begin{equation}
\left|
\frac{\phi^{[k]}(v,X)}{\phi(v,X)} 
\right| 
\leq C \tilde{v}^{\frac{(2m+1)k}{2m}}. 
\label{eqn:6.700}
\end{equation}
Moreover (\ref{eqn:6.550}),(\ref{eqn:6.70}) imply that 
$$ 
\left|
\frac{\partial^{\mu_0+1}}{\partial X^{\mu_0+1}} 
\frac{1}{\phi(v,X)}
\right|
\leq
C \tilde{v}^{\frac{(2m+1)(\mu_0+1)+m-1}{2m}} 
e^{-a(\tilde{X})\tilde{v}}. 
$$
Now we admit  the following lemma.
%%%%%%%%%%%%%%%%%%%%%%%%%%%%%%%%%%%%%%%%%%%%%%%%
\begin{lemma}\quad
$a(\tilde{X})\geq a(0)=a.$ 
\end{lemma}
%%%%%%%%%%%%%%%%%%%%%%%%%%%%%%%%%%%%%%%%%%%%%%%%%
We remark that the constant $a$ is as in Lemma 6.1. 
The above lemma implies 
\begin{equation}
\left|
\frac{\partial^{\mu_0+1}}{\partial X^{\mu_0+1}} 
\frac{1}{\phi(v,X)}
\right|
\leq
C \tilde{v}^{\frac{(2m+1)(\mu_0+1)+m-1}{2m}} 
e^{-a \tilde{v}}. 
\label{eqn:6.71}
\end{equation}
Finally   substituting (\ref{eqn:6.71}) into (\ref{eqn:6.54}), 
we can obtain the estimate of $r_{\mu_0}$ in Lemma 6.5. 
\qed  

%%%%%%%%%%%%%%%%%%%%%%%%%%%%%%%%%%%%%%%%%%%%%%%%

{\it Proof of Lemma 6.7}. \,\,\,
The definition of $\alpha(\tilde{X})$ in (\ref{eqn:6.58}) 
implies that 
\begin{eqnarray*} 
a'(\tilde{X})
\!\!\!&=&\!\!\!
-\alpha'(\tilde{X})\cdot \alpha(\tilde{X})^{2m-1}
[\hat{g}'(\tilde{X}a(\tilde{X}))\tilde{X}\alpha(\tilde{X}) 
+2m \hat{g}(\tilde{X}a(\tilde{X}))] \\
&& \quad \quad \quad +\alpha'(\tilde{X}) 
-g'(\tilde{X}\alpha(\tilde{X}))\alpha(\tilde{X})^{2m+1}\\
\!\!\!&=&\!\!\! 
-\hat{g}'(\tilde{X}\alpha(\tilde{X}))\alpha(\tilde{X})^{2m+1}.  
\end{eqnarray*} 

Since the condition $xg'(x)\leq 0$ in (\ref{eqn:2.1}) implies 
$\tilde{X}a'(\tilde{X})=
-\tilde{X} \hat{g}'(\tilde{X}\alpha(\tilde{X})) 
\cdot\alpha(\tilde{X})^{2m+1} 
\geq 0$, 
$a(\tilde{X})$ takes the minimum value when $\tilde{X}=0$. 
It is easy to check that $a(0)=a$. 
\qed

\subsection{Proof of Proposition 6.2}

(i)\,\,\,
Recall the definition of the function 
$K_{\mu}^{\langle 2 \rangle}$:  
$$
K_{\mu}^{\langle 2 \rangle}(\tau,\xi)
=
\int_{\xi}^{1}
e^{-s^{2m}}L_{\mu}(t_0 s)s^{4m+1-\mu}ds, 
$$
where 
$$
L_{\mu}(u) 
=
\int_{-\infty}^{\infty}
e^{uv}a_{\mu}(v) dv. 
$$ 
We remark that $L_{\mu}$  extends to an entire function. 

%\vspace{.5em}

By the residue formula, 
we have 
\begin{equation}
K_{\mu}^{\langle 2 \rangle}(\tau,\xi)
=
h_{\mu}(t_0) \xi^{-4m-2+\mu} \log \xi 
+ \tilde{h}_{\mu}(t_0,\xi),
%\label{eqn:6.74}
\end{equation}
where
\begin{equation}
h_{\mu}(t_0)=
\frac{1}{2 \pi i} 
\oint_{|\zeta|=\delta} e^{-\zeta^{2m}}L_{\mu}(t_0 \zeta)
\zeta^{4m+1-\mu} d\zeta 
\label{eqn:6.74}
\end{equation}
for $\delta$ is a small positive integer and 
$\tilde{h}_{\mu} \in C^{\infty}([0,1]\times[0, \varepsilon))$. 
Here (\ref{eqn:6.74}) implies that $h_{\mu}\equiv 0$ for 
$0 \leq \mu \leq 4m+1$ and 
$h_{\mu} \in C^{\infty}([0,1])$ for $\mu \geq 4m+2$. 
By Lemma 6.3, we can obtain (i) in Proposition 6.2. 

%\vspace{.5em}
%%%%%%%%%%%%%%%%%%%%%%%%%%%%%%%%
(ii)\,\,\,
%Recall the definition of the function 
%$\tilde{K}_{\mu_0}^{\langle 2 \rangle}$ and 
%$\tilde{L}_{\mu_0}$: 
%\begin{eqnarray}
%&&
%\tilde{K}_{\mu_0}^{\langle 2 \rangle}(\tau,\xi)
%=
%\int_{\xi}^{1}
%e^{-s^{2m}} \tilde{L}_{\mu_0}(t_0,\xi;s)s^{4m-\mu_0}ds, \\
%&&
%\tilde{L}_{\mu_0}(t_0,\xi;s)
%=
%\int_{-\infty}^{\infty}
%e^{t_0 sv}r_{\mu_0}(v,\xi s^{-1})dv. 
%\end{eqnarray}
Changing the integral variable, 
we have 
\begin{eqnarray*}
&&
\xi^{-4m-2+\mu_0} \tilde{K}_{\mu_0}^{\langle 2 \rangle}(\tau,\xi) 
=
\int_1^{\xi^{-1}} e^{-\xi^{2m}u^{2m}}
\tilde{L}_{\mu_0}(t_0,\xi,\xi u) u^{4m-\mu_0}du, \\ 
&&
\tilde{L}_{\mu_0}(t_0,\xi,\xi u)
=\int_{-\infty}^{\infty}  
e^{t_0\xi uv}r_{\mu_0}(v,u^{-1})dv. 
\end{eqnarray*}
Note that $r_{\mu_0}$ satisfies the inequality  in Lemma 6.5.  

%\vspace{.5em}

Keeping the above integrals in mind, 
we define the function $H$  by 
\begin{equation}
H(\alpha,\beta,\gamma,\delta) 
=t_0^{\delta}\xi^{\alpha} \int_1^{\xi^{-1}} 
e^{-\xi^{2m}u^{2m}}\frac{\partial^{\gamma}}{\partial X^{\gamma}} 
I(t_0 \xi u) u^{-\beta-1} du, 
\label{eqn:6.77}
\end{equation}
with 
\begin{equation}
I(X)=\int_{-\infty}^{\infty}e^{Xv} r(v)dv, 
\label{eqn:6.770}
\end{equation}
where 
$\alpha,\beta,\gamma(\geq 0),\delta$ are integers and 
the function $r$ satisfies 
$|r(v)| \leq C e^{-c|v|^{\frac{2m}{2m-1}}}$ 
for some positive constants $c$, $C$. 
Note that $H$ is a function of $(t_0,\xi)$.

By a direct computation, we have 
\begin{eqnarray}
\frac{\partial}{\partial t_0} 
H(\alpha,\beta,\gamma,\delta) 
\!\!\!&=&\!\!\!
H(\alpha+1,\beta-1,\gamma+1,\delta),\label{eqn:6.78}\\ 
\frac{\partial}{\partial \xi} 
H(\alpha,\beta,\gamma,\delta)
\!\!\!&=&\!\!\! 
-e^{-1}t_0^{\delta}\,\xi^{\alpha+\beta-1} 
\frac{\partial^{\gamma}}{\partial X^{\gamma}}I(\xi)\nonumber
+\alpha H(\alpha-1,\beta,\gamma,\delta)\\
&&\!\!\!
-2m H(\alpha+2m-1,\beta-2m,\gamma,\delta)\nonumber\\ 
&&\!\!\!
+H(\alpha,\beta-1,\gamma+1,\delta+1). 
\label{eqn:6.79}
\end{eqnarray} 
Since $\frac{\partial^{\gamma}}{\partial X^{\gamma}}I$ 
is bounded on $[0,1]$, we have 
\begin{equation}
|H(\alpha,\beta,\gamma,\delta)| 
\leq C|\xi|^{\alpha+\beta}, 
\label{eqn:6.80} 
\end{equation} 
by (\ref{eqn:6.77}).
By induction, we have 
$$
\left|
\frac{\partial^k}{\partial t_0^k}\frac{\partial^l}{\partial \xi^l}
H(\alpha,\beta,\gamma,\delta)
\right| 
\leq C|\xi|^{\alpha+\beta-l}, 
$$
by (\ref{eqn:6.78}),(\ref{eqn:6.79}),(\ref{eqn:6.80}). 
 
Now if we replace $r(v)$ by $r_{\mu_0}(v,Y)$ in (\ref{eqn:6.770}), 
then 
$\xi^{-4m-2+\mu_0} \tilde{K}_{\mu_0}^{\langle 2 \rangle}(\tau,\xi)=
H(0,-4m+\mu_0-1,0,0)$. 
Therefore if $\mu_0>4m+1+l$, then  
$\frac{\partial^k}{\partial t_0^k} 
 \frac{\partial^l}{\partial \xi^l}  
 K_{\mu_0}^{\langle 2 \rangle}(\tau,\xi)$ is a 
continuous function of $(t_0, \xi)\in [0,1]\times[0,\varepsilon)$.  
Therefore we can obtain (ii) in Proposition 6.2 by Lemma 6.3. 

This completes the proof of Proposition 6.2. 
\qed

%%%%%%%%%%%%%%%%%%%%%%%%%%%%%%%%%%%%%%%%%%%%%%%%%%%%%%%%%%%%%

\section{The Szeg\"o kernel} 

Let $\Omega_f$ be a tube domain 
satisfying the condition in Section 2. 
Let $H^2(\Omega_f)$ be the subspace of $L^2(\Omega_f)$ 
consisting of holomorphic functions $F$ on $\Omega_f$ 
such that 
$$
\sup_{\epsilon>0} 
\int_{\partial\Omega_f}
|F(z_1,z_2+i\epsilon)|^2  
d\sigma(z) <\infty, 
$$
where $d\sigma$ is the measure on $\partial\Omega_f$ given by 
Lebesgue measure on  ${\bf C}\times{\bf R}$  
when we identify $\partial\Omega_f$ with ${\bf C}\times{\bf R}$ 
($(z,t+if({\rm Im}z))\mapsto (z,t)$). 
The Szeg\"o projection is the orthogonal projection 
${\bf S} : L^2(\partial\Omega_f)\to H^2(\Omega_f)$ 
and we can write 
$$
{\bf S}F(z)=\int_{\partial\Omega_f} 
S(z,w)F(w)d\sigma(w), 
$$
where $S:\Omega_f\times\Omega_f\to {\bf C}$ is the 
{\it Szeg\"o kernel} of the domain $\Omega_f$. 
We are interested in the restriction of the Szeg\"o kernel 
on the diagonal, so we write $S(z)=S(z,z)$. 

%\vspace{.5em}

The Szeg\"o kernel of $\Omega_f$ has an integral representation : 
$$ 
S(z)=\frac{1}{(4\pi)^2} \int\!\!\!\int_{\Lambda^{\ast}} 
e^{-x\zeta_1-y\zeta_2} 
\frac{1}{D(\zeta_1,\zeta_2)} d\zeta_1d\zeta_2, 
$$
where $(x,y)=({\rm Im}z_1,{\rm Im}z_2)$ and 
$D(\zeta_1,\zeta_2)$ is as in Section 3 (\ref{eqn:3.4}). 

%\vspace{.4 em}

We also give an asymptotic expansion 
of the Szeg\"o kernel of $\Omega_f$. 
The theorem below can be obtained in a fashion similar to 
the case of the Bergman kernel, 
so we omit the proof.

%%%%%%%%%%%%%%%%%%%%%%%%%%%%%
%%%%%%%%%%%%%%%%%%%%%%%%%%%%%
%%%%%%%%%%%%%%%%%%%%%%%%%%%%%

\begin{thm}
The Szeg\"o kernel of $\Omega_f$ has the form in some  
neighborhood of  $z^0$: 
$$
S(z)=\frac{\Phi^S(\tau, \varrho^{\frac{1}{m}})}
{\varrho^{1+\frac{1}{m}}} 
+\tilde{\Phi}^S(\tau,\varrho^{\frac{1}{m}}) \log \varrho^{\frac{1}{m}},
$$
where 
$\Phi^S \in C^{\infty}((0,1]\times [0,\varepsilon))$ and 
$\tilde{\Phi}^S \in C^{\infty}([0,1]\times [0,\varepsilon))$,  
with some $\varepsilon>0$. 

Moreover 
$\Phi^S$ is written in the form on the set 
$\{\tau> \alpha \varrho^{\frac{1}{2m}}\}$ with 
some $\alpha>0$: 
for every nonnegative integer $\mu_0$  
$$
\Phi^S(\tau,\varrho^{\frac{1}{m}}) 
=\sum_{\mu=0}^{\mu_0} c_{\mu}^S(\tau) 
\varrho^{\frac{\mu}{m}} 
+ R_{\mu_0}^S(\tau, \varrho^{\frac{1}{m}})
\varrho^{\frac{\mu_0}{m}+\frac{1}{2m}}, 
$$
where 
$$
c_{\mu}^S(\tau)
= \frac{\varphi_{\mu}^S(\tau)}
{\tau^{2+2\mu}} 
+\psi_{\mu}^S(\tau)\log \tau, 
$$  
for 
$\varphi_{\mu}^S,\psi_{\mu}^S \in 
C^{\infty}([0,1])$,  $\varphi_0^S$ is positive on $[0,1]$ 
and $R_{\mu_0}^S$ satisfies 
$
|R_{\mu_0}^S(\tau,\varrho^{\frac{1}{m}})|
 \leq C_{\mu_0}^S[\tau-\alpha\varrho^{\frac{1}{2m}}]^{-3-2\mu_0}  
$ 
for some positive constant $C_{\mu_0}^S$. 
\end{thm}
%%%%%%%%%%%%%%%%%%%%%%%%%%
%%%%%%%%%%%%%%%%%%%%%%%%%%

%%%%%%%%%%%%%%%%%%%%%%%%%%%%%%%%%%%%%%%%%%%%%%%%%%%%%%%%%

%\end{large}

\end{document}